\documentclass[preprint,authoryear]{elsarticle}
\usepackage[latin1]{inputenc}
\usepackage{amsmath}
\usepackage{amsfonts}
\usepackage{amssymb}
\usepackage{amsthm}
\usepackage{graphicx}
\usepackage{placeins}
\usepackage{natbib}   

\def \X {\mathcal{X}}
\def \C {\mathcal{C}}
\def \RR {\mathbb{R}}
\def\EE {\mathbb{E}}
\def \N {\mathcal{N}}
\def \Z {\mathbb{Z}}

\newtheorem{prop}{Proposition}[section]

\newtheorem{lem}[prop]{Lemma}
\newtheorem{defi}[prop]{Definition}
\newtheorem*{lem*}{Lemma}

\newcommand{\argmin}{\mathop{\mathrm{argmin}}}

\title{Cross Validation and Maximum Likelihood estimations of hyper-parameters of Gaussian processes with model misspecification}
              
\tnotetext[t1]{A supplementary material is attached in the electronic version of the article.}

\author[fb1,fb2]{Fran\c cois Bachoc\corref{cor1}}

\cortext[cor1]{Corresponding author: Fran\c cois Bachoc \\
CEA-Saclay, DEN, DM2S, STMF, LGLS, F-91191 Gif-Sur-Yvette, France \\
Phone: +33 (0) 1 69 08 97 91 \\
Email: francois.bachoc@cea.fr}                                                

\address[fb1]{CEA-Saclay, DEN, DM2S, STMF, LGLS, F-91191 Gif-Sur-Yvette, France \\}
\address[fb2]{Laboratoire de Probabilit\'es et Mod\`eles Al\'eatoires, Universit\'e Paris VII \\
Site Chevaleret, case 7012, 75205 Paris cedex 13}

\begin{document}

\begin{abstract}
The Maximum Likelihood (ML) and Cross Validation (CV) methods for estimating covariance
hyper-parameters are compared, in the context of Kriging with a misspecified covariance structure.
A two-step approach is used.
First, the case of the estimation of a single variance hyper-parameter is addressed, for which the fixed correlation function is misspecified.
A predictive variance based quality
criterion is introduced and a closed-form expression of this criterion is derived. It is shown that when the correlation function is misspecified, the CV does better compared to
ML, while ML is optimal when the model is well-specified. In the second step, the results of the first step are extended to the case when the hyper-parameters of the
correlation function are also estimated from data. 
\end{abstract}


\maketitle

\newpage

\section{Introduction}

Kriging models \citep{TVRA,ISDSTK,GPML} consist in interpolating the values of a Gaussian random field given observations at a finite set of observation points.
They have become a popular method for a
large range of applications, such as numerical code approximation \citep{DACE,TDACE} and calibration \citep{CCMMO} or global optimization \citep{EGOEBBF}.

One of the main issues regarding Kriging is the choice
of the covariance function.
A Kriging model yields an unbiased predictor, with minimal variance
and a correct predictive variance, if the covariance function used in the model coincides with that of the random field from which the
observations stem.
Significant results concerning the influence of a misspecified covariance function on
Kriging predictions are available \citep{AEPRFMCF,BELPUICF,UAOLPRFUISOS}. These results are based on the fixed-domain asymptotics
\citep[p.62]{ISDSTK}, that is to say the case when the infinite sequence of observation points is dense in a bounded domain.
In this setting, the results in \citep{AEPRFMCF,BELPUICF,UAOLPRFUISOS} state that if the
true and assumed covariance functions yield
equivalent Gaussian measures (see e.g \citet[p.110]{ISDSTK} and \citet[p.63]{GRP}), then there is asymptotically no loss using the incorrect covariance function.
The asymptotic optimality holds for both the predictive means and variances.

Hence, in the fixed-domain asymptotics framework, it is sufficient to estimate a covariance function equivalent to the true one.
Usually, it is assumed that the covariance function belongs to a given parametric family. In this case, the estimation boils down to estimating the corresponding
parameters, that are called "hyper-parameters".
Because of the above theoretical
results, it is useful to distinguish between microergodic and non-microergodic hyper-parameters. Following the definition
in \citet{ISDSTK}, an hyper-parameter is microergodic if two covariance functions are orthogonal whenever they differ for it (as in \citet{ISDSTK}, we say that
two covariance functions are orthogonal if the two underlying Gaussian measures are orthogonal). Non-microergodic hyper-parameters
cannot be consistently estimated but have no asymptotic influence on Kriging predictions, as shown in \citet{IEAEIMBG} in the Mat\'ern case. There is a
fair amount of literature on the consistent estimation of microergodic hyper-parameters, particularly using the Maximum Likelihood (ML) method.
Concerning the isotropic Mat\'ern model, with
fixed regularity parameter $\nu$ and free variance $\sigma^2$ and correlation length $\ell$, the microergodic ratio $\frac{\sigma^2}{\ell^{2\nu}}$
can be consistently estimated by ML for dimension $d \leq 3$ \citep{IEAEIMBG},
with asymptotic normality for $d=1$ \citep{FDAPTMLE}. Both $\sigma^2$ and $\ell$ can be consistently estimated for $d>4$ \citep{CSSVGRF}.
For the multiplicative Mat\'ern
model with $\nu = \frac{3}{2}$, both $\sigma^2$ and the $d$ correlation length hyper-parameters $\ell_1,...,\ell_d$ are consistently estimated by ML for $d \geq 3$
\citep{FDASMTGRF}. For the Gaussian
model, the $d$ correlation lengths are consistently estimated by ML \citep{ESCMSGRFM}. Finally for the multiplicative exponential model,
the microergodic ratio $\frac{\sigma^2}{\ell}$ is consistently estimated by ML with asymptotic normality for $d=1$ \citep{APMLEDGP}.
All hyper-parameters $\sigma^2,\ell_1,...,\ell_d$ are consistently estimated by ML with asymptotic normality for $d > 1$ \citep{MLEPUSSS}.

We believe that the fixed-domain asymptotics does not solve completely the issue of the estimation of the covariance function.
The first point is that the above theoretical results
are asymptotic, while one could be interested in finite-sample studies, such as the numerical experiments performed in \citet[ch.6.9]{ISDSTK},
and the detailed study of the exponential covariance in \citet{TRAFSS}.
The second point
is that one may not be able to asymptotically estimate a covariance function equivalent to the true one.
Indeed, for instance, for two covariance functions of the isotropic Mat\'ern class to be equivalent,
it is necessary that their regularity parameters are equal (see the one-dimensional
example in \citet[p.136]{ISDSTK}). Yet, it is common practice, especially for the analysis of computer experiment data, to enforce the regularity parameter
to an arbitrary value (see e.g \citet{UKMADCM}). The interplay between the misspecification of the regularity parameter
and the prediction mean square error is not trivial. The numerical experiments in \citet[ch.5.3.3]{MCSLMMNA}
show that a misspecified regularity parameter can have dramatic consequences.

The elements pointed out above justify addressing the
case when a parametric estimation is carried out, within a covariance function set, and when the true underlying covariance
function does not belong to this set. We call this the model misspecification case. In this context, we study
the Cross Validation (CV) estimation method \citep{PACHGP,KCVMSD}, and compare it with ML.
This comparison has been an area of active research. Concerning theoretical
results, \citet{CGCVMMLEPSP} showed that for the estimation of a signal-to-noise ratio parameter of a Brownian motion, CV has twice the asymptotic variance of ML
in a well-specified case.
Several numerical results are also available, coming either from Monte Carlo studies as in \citet[ch.3]{TDACE} or deterministic studies as in \citet{UKMADCM}.
In both the above studies, the interpolated functions are smooth, and the covariance structures are adapted, being Gaussian in \citet{UKMADCM} and having a
free smoothness parameter in \citet{TDACE}.

We use a two-step approach.
In the first step, we consider a parametric family of stationary covariance functions in which only the global variance hyper-parameter
is free. In this framework, we carry out a detailed and quantitative finite-sample comparison, using the closed-form
expressions for the estimated variances for both the ML and CV methods. For the second step we study the
general case in which the global variance hyper-parameter and the correlation hyper-parameters are free and estimated from data.
We perform extensive numerical experiments on analytical functions, with various misspecifications, and we compare
the Kriging models obtained with the ML and CV estimated hyper-parameters.

The paper is organized as follows. In Section \ref{section: framework_variance}, we detail the statistical framework for the estimation of a single variance
hyper-parameter, we introduce an original
quality criterion for a variance estimator, and we give a closed-form formula of this criterion for a large family of estimators. In Section
\ref{section: LOOML} we introduce the ML and Leave-One-Out (LOO) estimators of the variance hyper-parameter. In
Section \ref{section: numerical_results_variance} we numerically apply the closed-form formulas of Section \ref{section: framework_variance} and we study their dependences
with respect to
model misspecification and number of observation points. We highlight our main result that when the correlation model is
misspecified, the CV does better compared to ML. Finally in Section \ref{section: analyticalfunctions} we illustrate this result on the
Ishigami analytical function and then generalize it, on the Ishigami and Morris analytical functions, to the case where the correlation hyper-parameters are estimated as well.

\section{Framework for the variance hyper-parameter estimation} \label{section: framework_variance}

We consider a Gaussian process $Y$, indexed by a set $\X$. $Y$ is stationary, with unit variance, and its correlation function is denoted by $R_1$.
A Kriging model is built for $Y$, for which it is assumed that $Y$ is centered and that its covariance function belongs to the set $\C$, with
\begin{equation} \label{eq: family_covariance}
\C = \left\{ \sigma^2 R_2, \sigma^2 \in \RR^+ \right\},
\end{equation}
with $R_2(x)$ a given stationary correlation function. Throughout this paper, $\EE_i$, $var_i$, $cov_i$ and $\sim_i$, $i \in \{1,2\}$, denote means, variances, covariances
and probability laws taken with respect to the distribution of $Y$ with mean zero, variance one, and the correlation function $R_i$.
We observe $Y$ on the points $x_1,...,x_n \in \X$. 
In this framework, the hyper-parameter $\sigma^2$ is estimated from the data $y = (y_1,...,y_n)^t = (Y(x_1),...,Y(x_n))^t$ using an
estimator $\hat{\sigma}^2$.
This estimation does not affect the Kriging prediction of $y_0 = Y(x_0)$, for a new point $x_0$, which
we denote by $\hat{y}_0$:
\begin{equation} \label{eq: hat_y0}
\hat{y}_0 := \EE_2 ( y_0 | y ) =  \gamma_2^t \mathbf{\Gamma}_2^{-1} y,
\end{equation}
where $(\gamma_i)_j = R_i( x_j - x_0)$ and $(\mathbf{\Gamma}_i)_{j,k} = R_i(x_j - x_k)$, $i \in \{1,2\}$, $1 \leq  j,k \leq n $.
The conditional mean square error of this non-optimal prediction is, after
a simple calculation,
\begin{equation} \label{eq: riskdepred}
\EE_1 \left[ ( \hat{y}_0 - y_0 )^2 | y \right] = ( \gamma_1^t \mathbf{\Gamma}_1^{-1} y - \gamma_2^t \mathbf{\Gamma}_2^{-1} y )^2 + 1 - \gamma_1^t \mathbf{\Gamma}_1^{-1} \gamma_1.
\end{equation}
However, using the covariance family $\C$, we use the classical Kriging predictive variance
expression $\hat{\sigma}^2  c_{x_0}^2$, with
\begin{equation} \label{eq: cx0}
c_{x_0}^2  := var_2( y_0 | y ) = 1 - \gamma_2^t \mathbf{\Gamma}_2^{-1} \gamma_2.
\end{equation}
As we are interested in the accuracy of the predictive variances obtained from an estimator $\hat{\sigma}^2$, the following notion of Risk can be formulated. 
\begin{defi} \label{def: risk}
For an estimator $\hat{\sigma}^2$ of $\sigma^2$, we call Risk at $x_0$ and denote by $\mathcal{R}_{ \hat{\sigma}^2 , x_0 }$ the quantity
\[
\mathcal{R}_{ \hat{\sigma}^2 , x_0 } = \EE_1 \left[  \left(  \EE_1 \left[ ( \hat{y}_0 - y_0 )^2 | y \right] -  \hat{\sigma}^2  c_{x_0}^2 \right)^2  \right].
\]
\end{defi}
If $\mathcal{R}_{ \hat{\sigma}^2 , x_0 }$ is small, then this means that the predictive variance $\hat{\sigma}^2 c_{x_0}^2$ is a correct prediction of the
conditional mean square error
\eqref{eq: riskdepred} of the
prediction $\hat{y}_0$. Note that when $R_1 = R_2$ the minimizer of the Risk at every $x_0$ is $\hat{\sigma}^2 =1$.
When $R_1 \neq R_2$, an estimate of $\sigma^2$ different from $1$ can improve the predictive variance, partly compensating for the
correlation function error. \\
To complete this section, we give the closed-form expression of the Risk of an estimator that can be written as a quadratic form of the observations,
which is the case for all classical estimators, including the ML and CV estimators of Section \ref{section: LOOML}.

\begin{prop} \label{prop: main}
Let $\hat{\sigma}^2$ be an estimator of $\sigma^2$ of the form $y^t \mathbf{M} y$ with $\mathbf{M}$ an $n \times n$ matrix. Denoting
$f(\mathbf{A},\mathbf{B}) = tr(\mathbf{A})tr(\mathbf{B}) + 2 tr(\mathbf{A}\mathbf{B})$, for $\mathbf{A}$, $\mathbf{B}$ $n \times n$ real matrices,
$\mathbf{M}_{0} =  (\mathbf{\Gamma}_2^{-1} \gamma_2 - \mathbf{\Gamma}_1^{-1} \gamma_1  )( \gamma_2^t \mathbf{\Gamma}_2^{-1} -  \gamma_1^t \mathbf{\Gamma}_1^{-1} ) \mathbf{\Gamma}_1$,
$\mathbf{M}_1 =  \mathbf{M} \mathbf{\Gamma}_1$,
$c_1 = 1 - \gamma_1^t \mathbf{\Gamma}_1^{-1} \gamma_1$
and
$c_2 = 1 - \gamma_2^t \mathbf{\Gamma}_2^{-1} \gamma_2$,
we have:
\begin{eqnarray*}
 \mathcal{R}_{ \hat{\sigma}^2 , x_0 } & = & f(\mathbf{M}_{0},\mathbf{M}_{0}) + 2 c_1 tr(\mathbf{M}_{0}) - 2 c_2 f(\mathbf{M}_{0},\mathbf{M}_1) \\
 &  & + c_1^2 - 2 c_1 c_2 tr(\mathbf{M}_1) + c_2^2 f(\mathbf{M}_1,\mathbf{M}_1).
\end{eqnarray*}
\end{prop}

It seems difficult at first sight to conclude
from Proposition \ref{prop: main} whether one estimator is better than another given a correlation function error and a
set of observation points. Therefore,
in Section \ref{section: numerical_results_variance}, we numerically analyze the Risk for the ML and CV estimators of the variance
for several designs of experiments. Before that, we introduce the ML and CV estimators of $\sigma^2$.

\section{ML and CV estimation of the variance hyper-parameter} \label{section: LOOML}

In the framework of Section \ref{section: framework_variance}, the ML estimator $\hat{\sigma}^2_{ML}$ of $\sigma^2$ (see e.g \citet[p.66]{TDACE}) is
\begin{equation} \label{eq: hatsigmaML}
\hat{\sigma}_{ML}^2 = \frac{1}{n} y^t \mathbf{\Gamma}_2^{-1} y.
\end{equation}

Let us now consider the CV estimator of $\sigma^2$. The principle is that, given a value $\sigma^2$ specifying the covariance
function used among the set $\C$, we can, for $1 \leq i \leq n$, compute $\hat{y}_{i,-i} := \EE_2(y_i|y_1,...,y_{i-1},y_{i+1},...,y_n)$
and $\sigma^2 c_{i,-i}^2 := \sigma^2 var_2(y_i|y_1,...,y_{i-1},y_{i+1},...,y_n)$. The Cross Validation estimate of $\sigma^2$ is based
on the criterion 
\begin{equation} \label{eq: critLOO}
C_{LOO} = \frac{1}{n} \sum_{i=1}^n \frac{(y_i - \hat{y}_{i,-i} )^2}{\sigma^2 c_{i,-i}^2}.
\end{equation}
It is noted in \citet[p.102]{SSD} that if $\sigma^2 R_2$ is a correct estimate of the true covariance function,
then we should expect this criterion to be close to $1$.
In \ref{subsection: preuveCritLOO}, we show that if the observation points expand in a regular way, then for most classical correlation functions (the Mat\'ern family for example),
if $\sigma^2 R_2$ is the true covariance function, then \eqref{eq: critLOO} converges toward $1$ in the mean square sense.

Based on this rigorous result about the case of an infinite regular grid, we can conjecture that \eqref{eq: critLOO} should be close to $1$ in a general context
when the value of $\sigma^2$ is correct. Hence the idea of the CV estimation of $\sigma^2$ is to seek the value of $\sigma^2$ so that
this criterion is equal to $1$, which yields the estimator
\begin{equation} \label{eq: hatsigmaCV}
\hat{\sigma}_{CV}^2 = \frac{1}{n} \sum_{i=1}^n \frac{(y_i - \hat{y}_{i,-i} )^2}{ c_{i,-i}^2}.
\end{equation}

By means of the well-known virtual LOO formulas of Proposition \ref{prop: dubrule} (see e.g \citet[ch.5.2]{SS}), we obtain the following vector-matrix
closed-form expression of \eqref{eq: hatsigmaCV},
\[
\hat{\sigma}_{CV}^2 = \frac{1}{n} y^t \mathbf{\Gamma}_2^{-1} \left[ diag( \mathbf{\Gamma}_2^{-1} ) \right]^{-1} \mathbf{\Gamma}_2^{-1} y,
\]
with $diag( \mathbf{\Gamma}_2^{-1} )$ the matrix obtained by setting to $0$ all
non-diagonal terms of $\mathbf{\Gamma}_2^{-1}$.
In the supplementary material, we give a short reminder, with proofs, of the virtual LOO formulas given in proposition \ref{prop: dubrule}. 

\begin{prop} \label{prop: dubrule}
For $1 \leq i \leq n$:
\[
c_{i,-i}^2 = \frac{1}{(\mathbf{\Gamma}_2^{-1})_{i,i}}
\]
and
\[
y_i - \hat{y}_{i,-i} = \frac{1}{(\mathbf{\Gamma}_2^{-1})_{i,i}} ( \mathbf{\Gamma}_2^{-1} y )_i.
\]
\end{prop}

When $R_1 = R_2$, we see that ML is more efficient than CV. Indeed
\begin{equation} \label{eq: RiskR1=R2}
\mathcal{R}_{ \hat{\sigma}^2 , x_0 } =
\EE_1 \left( ( 1-\gamma_1^t \mathbf{\Gamma}_1^{-1}\gamma_1 ) - \hat{\sigma}^2  ( 1-\gamma_1^t \mathbf{\Gamma}_1^{-1}\gamma_1 ) \right)^2
= ( 1-\gamma_1^t \mathbf{\Gamma}_1^{-1}\gamma_1 )^2 \EE_1( ( \hat{\sigma}^2 -1 )^2 ),
\end{equation}
so that the Risk of definition \ref{def: risk} is proportional to the quadratic error in estimating the true $\sigma^2=1$.
We calculate $\EE_1( \hat{\sigma}^2_{ML} ) = \EE_1( \hat{\sigma}^2_{CV} ) = 1$, hence both
estimators are unbiased. Concerning their variances, on the one hand we recall in \ref{subsection: R1=R2} that $var_1( \hat{\sigma}^2_{ML} ) = \frac{2}{n}$, the 
Cram\'er-Rao bound for the estimation of $\sigma^2$, in the case when the true $\sigma^2$ is $1$
and the model for $Y$ is $\C$. On the other hand:
\begin{eqnarray*}
var_1( \hat{\sigma}^2_{CV} ) & = & \frac{2}{n^2} tr( \mathbf{\Gamma}_1^{-1} \left[ diag( \mathbf{\Gamma}_1^{-1} ) \right]^{-1} \mathbf{\Gamma}_1^{-1} \left[ diag( \mathbf{\Gamma}_1^{-1} ) \right]^{-1} ) \\
& = & \frac{2}{n^2} \sum_{i=1}^n \sum_{j=1}^n \frac{ ( \mathbf{\Gamma}_1^{-1} )_{i,j}^2 }{( \mathbf{\Gamma}_1^{-1} )_{i,i} ( \mathbf{\Gamma}_1^{-1} )_{j,j}  }.
\end{eqnarray*}
Hence
$var_1( \hat{\sigma}^2_{CV} ) \geq \frac{2}{n^2} \sum_{i=1}^n  \frac{ ( \mathbf{\Gamma}_1^{-1} )_{i,i}^2 }{( \mathbf{\Gamma}_1^{-1} )_{i,i} ( \mathbf{\Gamma}_1^{-1} )_{i,i}  } = \frac{2}{n}$,
the Cram\'er-Rao bound.

However $var_1( \hat{\sigma}^2_{CV} )$ is only upper-bounded by $2$ (because $\mathbf{\Gamma}_1^{-1}$ is a
covariance matrix). Furthermore $var_1( \hat{\sigma}^2_{CV} )$ can be arbitrarily close to $2$.
Indeed, with
\[
\mathbf{\Gamma}_1 = \mathbf{\Gamma}_2 =  \frac{n-1 + \epsilon}{n-1} \mathbf{I} - \frac{\epsilon}{n-1} \mathbf{J},
\]
where $\mathbf{J}$ is the $n \times n$ matrix with all coefficients being $1$, we obtain
\[
var_1( \hat{\sigma}^2_{CV} ) = \frac{2}{n} + \frac{2 (n-1)}{n} \frac{\epsilon^2}{ ( \epsilon + (n-1)(1-\epsilon) )^2 }.
\]
Hence, when $R_1 = R_2$, ML is more efficient to estimate the variance parameter. The object of the next section is to study the case $R_1 \neq R_2$ numerically.

\section{Numerical results for the variance hyper-parameter estimation} \label{section: numerical_results_variance}

All the numerical
experiments are carried out with the numerical software Scilab \citep{ESCS}. We use the Mersenne Twister pseudo
random number generator of M. Matsumoto and T. Nishimura, which is
the default pseudo random number generator in Scilab for large-size random simulations.

\subsection{Criteria for comparison} \label{subsection: criteria}

\paragraph{Pointwise criteria} 

We define two quantitative criteria that will be used to compare the ML and CV assessments of the predictive variance at prediction point $x_0$.

The first criterion is the Risk on Target Ratio (RTR),
\begin{equation} \label{eq: RTR}
RTR (x_0) = \frac{ \sqrt{ \mathcal{R}_{ \hat{\sigma}^2 , x_0 } } }{ \EE_1 \left[ ( \hat{y}_0 - y_0 )^2 \right]  },
\end{equation}
with $\hat{\sigma}^2$ being either $\hat{\sigma}^2_{ML}$ or $\hat{\sigma}^2_{CV}$. 

From Definition \ref{def: risk} we obtain
\begin{equation} \label{eq: RTR_expli}
RTR (x_0) = \frac{ \sqrt{   \EE_1 \left[  \left(  \EE_1 \left[ ( \hat{y}_0 - y_0 )^2 | y \right] -  \hat{\sigma}^2  c_{x_0}^2 \right)^2  \right]   } }{ \EE_1 \left[ ( \hat{y}_0 - y_0 )^2 \right]  }.
\end{equation}
The numerator of \eqref{eq: RTR_expli} is the mean square error in predicting the random quantity $\EE_1 \left[ ( \hat{y}_0 - y_0 )^2 | y \right]$
(the target in the RTR acronym) with
the predictor $ \hat{\sigma}^2  c_{x_0}^2$. The denominator of $\eqref{eq: RTR_expli}$ is, by the law of total expectation, the mean of the predictand
$\EE_1 \left[ ( \hat{y}_0 - y_0 )^2 | y \right]$. Hence, the RTR in \eqref{eq: RTR} is a relative prediction error, which is easily interpreted. 

We have the following bias-variance decomposition of the Risk, 
\begin{equation} \label{eq: risk_biais_variance}
\mathcal{R}_{ \hat{\sigma}^2 , x_0 }  = \left( \underbrace{   \EE_1 \left[ ( \hat{y}_0 - y_0 )^2 \right]  -  \EE_1 \left[ \hat{\sigma}^2 c_{x_0}^2 \right]   }_{\mbox{bias}}  \right)^2 +
\underbrace{ var_1 \left( \EE_1 \left[ ( \hat{y}_0 - y_0 )^2 | y  \right]  - \hat{\sigma}^2 c_{x_0}^2  \right) }_{\mbox{variance}}.
\end{equation}
Hence the second criterion is the Bias on Target Ratio (BTR) and is the relative bias
\begin{equation} \label{eq: BTR}
BTR(x_0) = \frac{ |   \EE_1 \left[ ( \hat{y}_0 - y_0 )^2 \right] - \EE_1 \left( \hat{\sigma}^2 c_{x_0}^2 \right)   | }{ \EE_1 \left[ ( \hat{y}_0 - y_0 )^2 \right] }.
\end{equation}
The following equation summarizes the link between RTR and BTR:
\begin{equation} \label{eq: link_RTR_BTR}
\left( \underbrace{RTR}_{\mbox{relative error}} \right)^2 = \left( \underbrace{BTR}_{\mbox{relative bias}} \right)^2  + \underbrace{ \frac{ var_1 \left( \EE_1 \left[ ( \hat{y}_0 - y_0 )^2 | y  \right]  - \hat{\sigma}^2 c_{x_0}^2  \right) }{ \EE_1 \left[ ( \hat{y}_0 - y_0 )^2 \right]^2 } }_{\mbox{relative variance}}  . 
\end{equation}
\paragraph{Pointwise criteria when $R_1=R_2$}

When $R_1=R_2$, $\EE_1 \left[ ( \hat{y}_0 - y_0 )^2 | y \right]$ does not depend on $y$. Therefore, the RTR and BTR simplify into
$RTR(x_0) = \sqrt{ \EE_1 \left[ ( \hat{\sigma}^2 - 1 )^2 \right] }$ and $BTR(x_0) = | 1 - \EE_1( \hat{\sigma}^2) |$.
Hence, the RTR and BTR are the mean square error and the bias in the estimation of the true variance $\sigma^2=1$, and $RTR^2 = BTR^2 + var_1( \hat{\sigma}^2 )$.

\paragraph{Integrated criteria} 

In this paragraph, we define the two integrated versions
of RTR and BTR over the prediction space $\X$. Assume $\X$ is equipped with a probability measure $\mu$. Then we define
\begin{equation}  \label{eq: IRTR}
IRTR = \sqrt{ \int_{\X} RTR^2(x_0) d \mu (x_0) }
\end{equation}  
and
\begin{equation}  \label{eq: IBTR}
IBTR = \sqrt{ \int_{\X} BTR^2(x_0) d \mu (x_0) }.
\end{equation}  
Hence we have the equivalent of \eqref{eq: link_RTR_BTR} for IRTR and IBTR: 
\begin{equation} \label{eq: link_IRTR_IBTR}
IRTR^2  =  IBTR^2 + 
\int_{\X} \frac{ var_1 \left( \EE_1 \left[ ( \hat{y}_0 - y_0 )^2 | y  \right]  - \hat{\sigma}^2 c_{x_0}^2  \right) }{ \EE_1 \left[ ( \hat{y}_0 - y_0 )^2 \right]^2 } d \mu(x_0). 
\end{equation}

\subsection{Designs of experiments studied} \label{subsection: doe}

We consider three different kinds of Designs Of Experiments (DOEs) of $n$ observation points
on the prediction space $\X = [0,1]^d$. 

The first DOE is the Simple Random Sampling (SRS) design and consists of $n$ independent
observation points with uniform distributions
on $[0,1]^d$. This design may not be optimal for a Kriging prediction point of view, as it is likely to contain relatively large areas without observation points.
However it is a convenient design for the estimation of covariance hyper-parameters because it may contain some points with small spacing. It is noted in
\citet[ch.6.9]{ISDSTK} that such points can dramatically improve the estimation of the covariance hyper-parameters. 

The second DOE is the Latin Hypercube Sampling Maximin (LHS-Maximin) design (see e.g \citet{TDACE}). This design is one of the
most widespread non-iterative designs in Kriging. To generate a LHS-Maximin design, we generate $1000$ LHS designs, and keep the one that maximizes
the min criterion (i.e. the
minimum distance between two different observation points). Let us notice that this is the method used by the Matlab function \textit{lhsdesign(...,'maximin',k)} which
generates $k$ LHS designs with default $k=5$.

The third DOE is a deterministic sparse regular grid. It is built according to the \citet{QIFTPCCF} sparse tensorization construction of the one-dimensional
regular grid $G = \{ \frac{1}{2^l} ,..., \frac{2^l - 1}{2^l} \}$ of level $l$. 

The three DOEs are representative of the classical DOEs that can be used for interpolation of functions, going from the most irregular ones
(SRS) to the most regular ones (sparse grid).

\subsection{Families of correlation functions studied} \label{subsection: families_cov}

We consider two classical families of stationary correlation functions. 
\begin{itemize}
\item The power-exponential correlation function family, parameterized by the vector of correlation lengths
$\ell = (\ell_1,...,\ell_d)$ and the power $p$. $R$ is power-exponential $(\ell,p)$ when
\begin{equation} \label{eq: pwrexp}
R(h) = \exp\left( - \sum_{i=1}^d \left( \frac{|h_i|}{\ell_i} \right)^{p} \right).
\end{equation}
\item The Mat\'ern correlation function family, parameterized by the vector of correlation lengths
$\ell = (\ell_1,...,\ell_d)$ and the regularity parameter $\nu$. $R$ is Mat\'ern $(\ell,\nu)$ when
\begin{equation} \label{eq: Rmat}
R(h) = \frac{1}{\Gamma(\nu)2^{\nu-1}} \left(  2 \sqrt{\nu} |h|_{\ell} \right)^{\nu} K_{\nu} \left( 2 \sqrt{\nu} |h|_{\ell} \right),
\end{equation}
with $ |h|_{\ell} = \sqrt{ \sum_{i=1}^d \frac{ h_i^2 }{\ell_i^2}  }$, $\Gamma$ the Gamma function and $K_{\nu}$ the modified Bessel function of second order.
See e.g \citet[p.31]{ISDSTK} for a presentation of the Mat\'ern correlation function.
\end{itemize}

\subsection{Influence of the model error} \label{subsection: inf_erreur_modele}
We study the influence of the model error, i.e. the difference between $R_1$ and $R_2$. For different pairs $R_1,R_2$,
we generate $n_p=50$ SRS and LHS learning samples, and the deterministic sparse grid (see Section \ref{subsection: doe}). We
compare the empirical means of the two integrated criteria IRTR and
IBTR (see Section \ref{subsection: criteria}) for the different
DOEs and for ML and CV.
IRTR and IBTR are calculated on a large test sample of size $5000$. We take $n=70$ for the learning sample size (actually $n=71$ for the regular grid)
and $d=5$ for the dimension.

For the pairs $R_1,R_2$, we consider the three following cases. First, $R_1$ is power-exponential $((1.2,...,1.2),1.5)$ and
$R_2$ is power-exponential $((1.2,...,1.2),p_2)$ with varying $p_2$.
Second, $R_1$ is Mat\'ern $((1.2,...,1.2),1.5)$ and $R_2$ is Mat\'ern $((1.2,...,1.2),\nu_2)$ with varying $\nu_2$.
Finally, $R_1$ is Mat\'ern $((1.2,...,1.2),1.5)$ and $R_2$ is Mat\'ern $((\ell_2,...,\ell_2),1.5)$ with varying $\ell_2$.

\begin{figure}[]
\centering
 \hspace*{-2cm}
\begin{tabular}{c c c}
\includegraphics[width=5cm,angle=0]{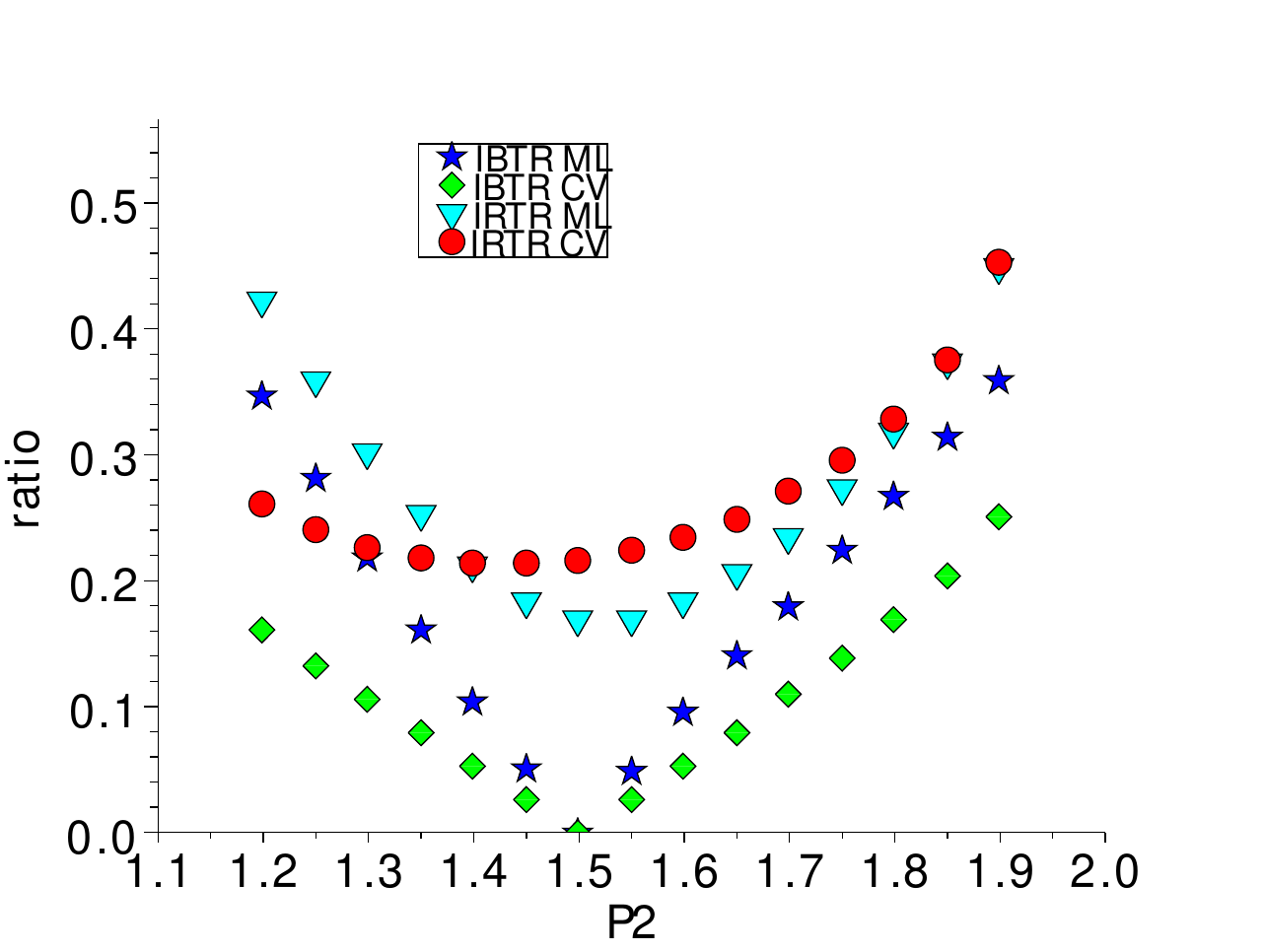} & \includegraphics[width=5cm,angle=0]{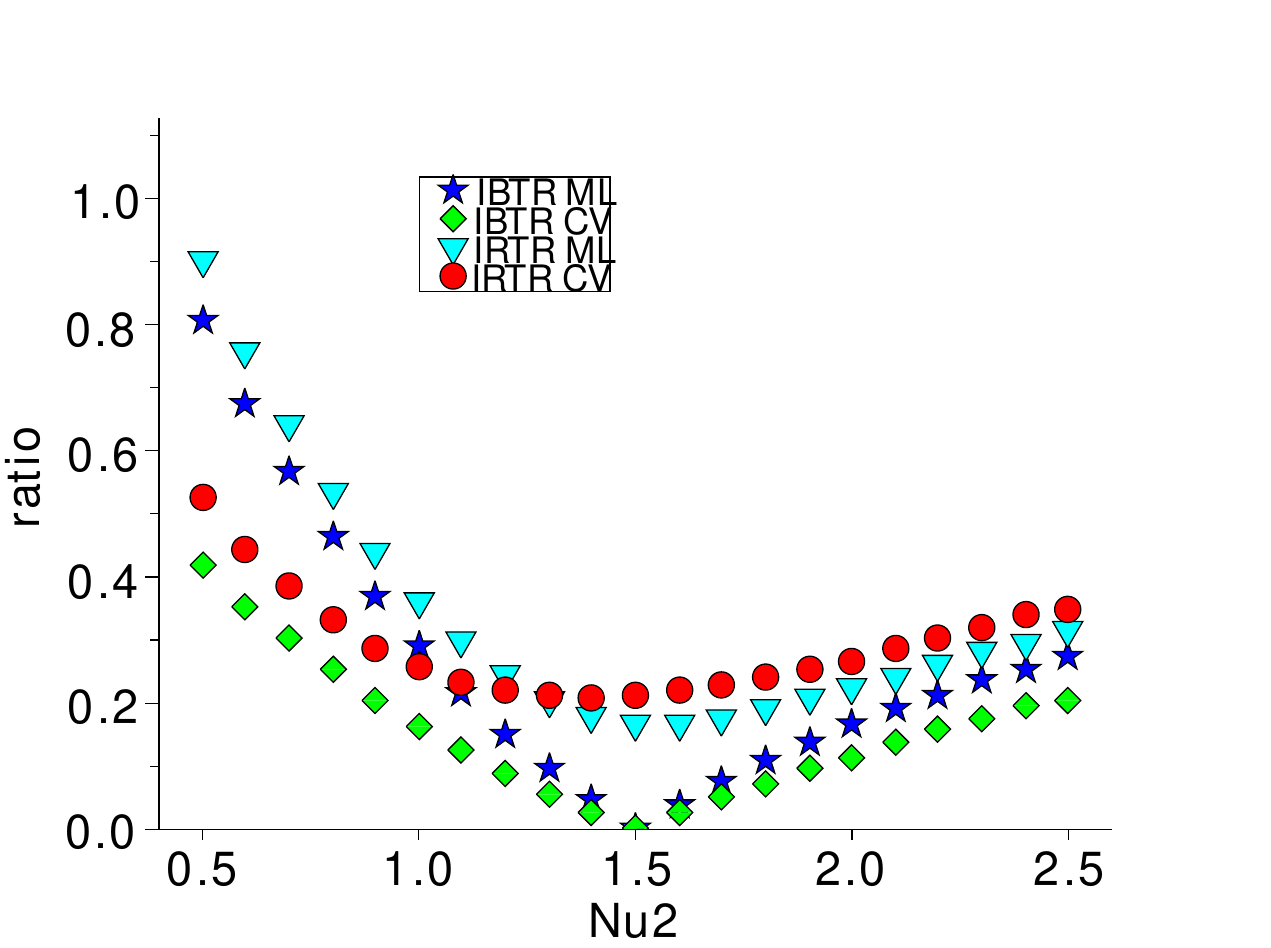}
& \includegraphics[width=5cm,angle=0]{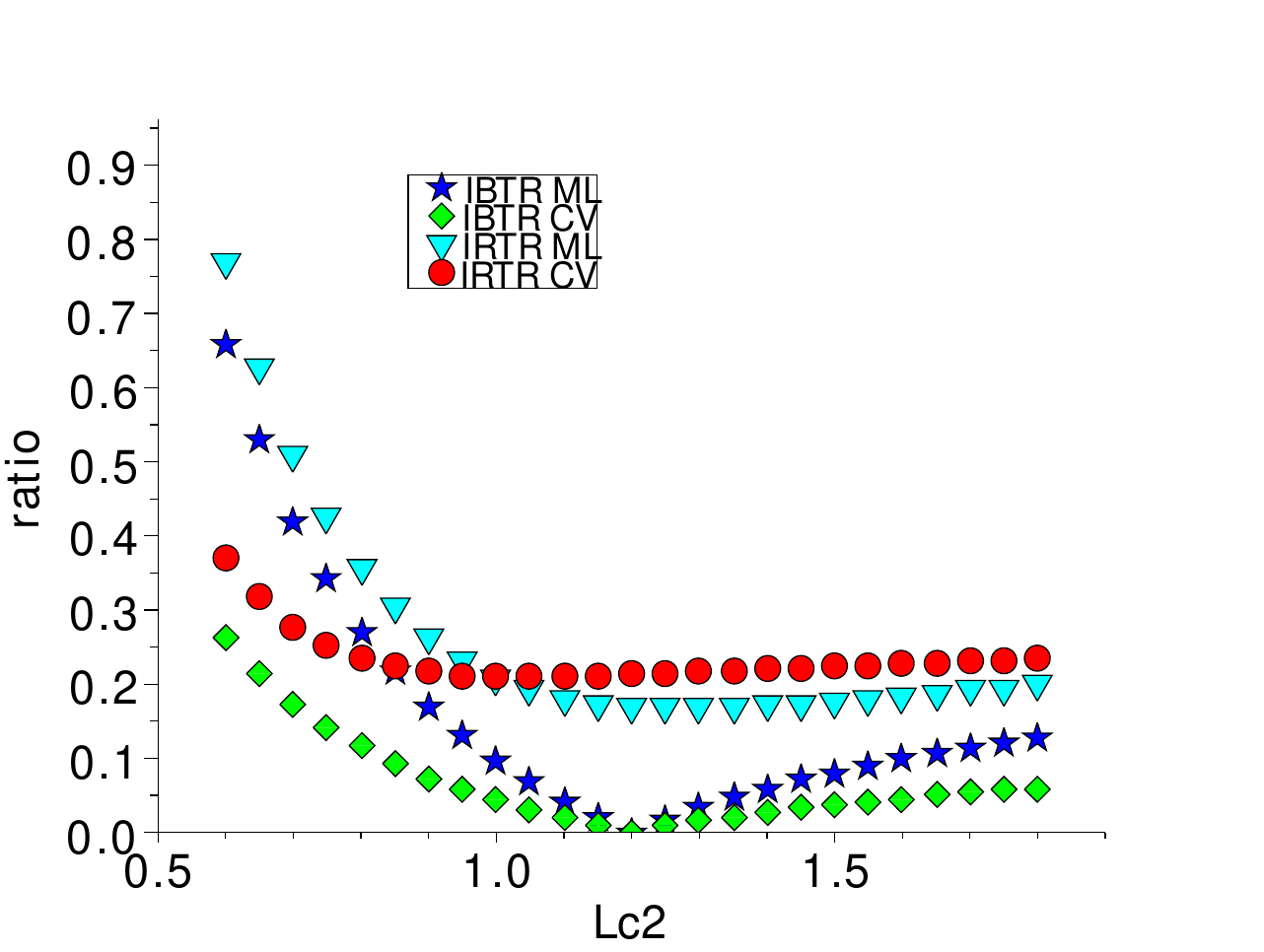}
\end{tabular}
\caption{Influence of the model error for the SRS DOE (see Section \ref{subsection: doe}). Plot of the IRTR and IBTR integrated criteria (see Section
\ref{subsection: criteria}) for ML and CV. Left: power-exponential correlation function with error on the exponent, the true exponent is $p_1=1.5$ and the model exponent
$p_2$
varies in $[1.2,1.9]$. Middle: Mat\'ern correlation function with error on the regularity parameter, the true regularity parameter is $\nu_1=1.5$ and the model
regularity parameter $\nu_2$ varies in $[0.5,2.5]$. Right: Mat\'ern correlation function with error on the correlation length, the true correlation length is $\ell_1=1.2$ and the
model correlation length $\ell_2$ varies in $[0.6,1.8]$. ML is optimal when there is no model error while CV is more robust to model misspecifications.}
\label{fig: IntegreSRS}
\end{figure}

On Figure \ref{fig: IntegreSRS}, we plot the results for the SRS DOE. We clearly see that when the model error becomes large,
CV becomes more efficient than ML in the sense of IRTR. Looking at
\eqref{eq: link_IRTR_IBTR}, one can see that the IRTR is composed of IBTR and of an integrated relative variance term. When $R_2$ becomes different from $R_1$, the IBTR
contribution increases faster than the integrated relative variance contribution,
especially for ML. Hence, the main reason why CV is more robust than ML to model misspecification is that its bias increases more slowly with the model
misspecification.

\begin{figure}[]
\centering
 \hspace*{-2cm}
\begin{tabular}{c c c}
\includegraphics[width=5cm,angle=0]{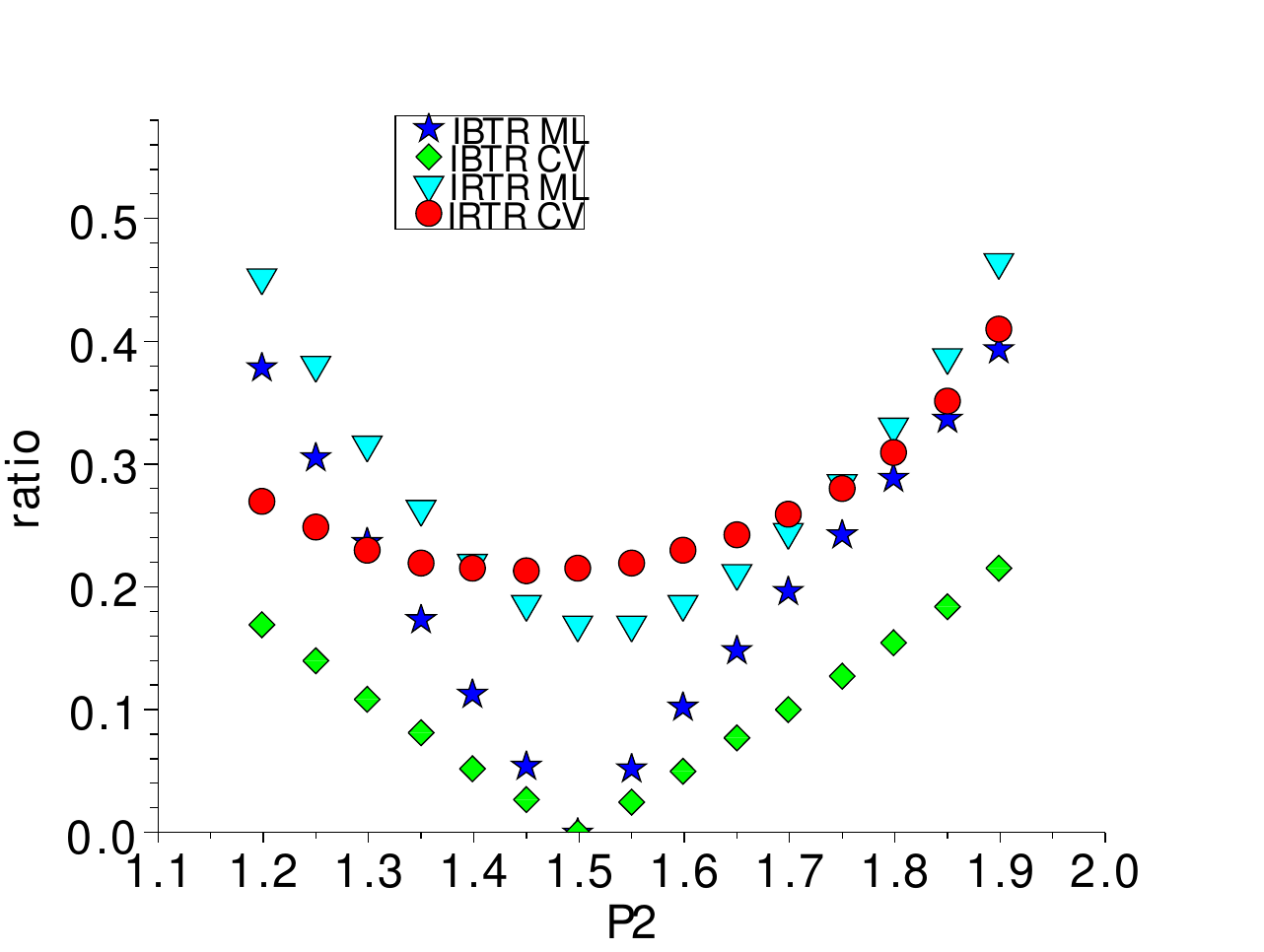} & \includegraphics[width=5cm,angle=0]{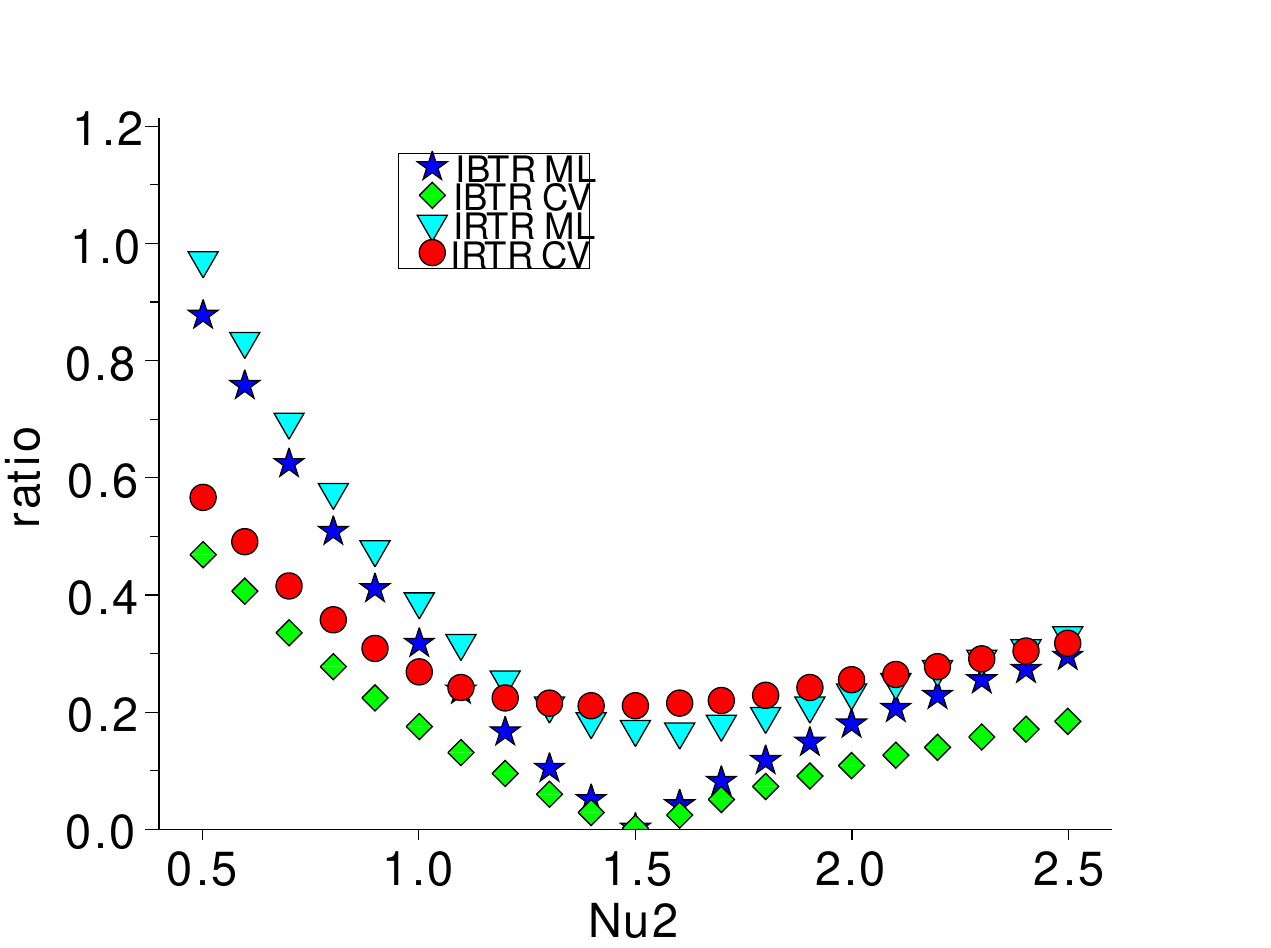}
& \includegraphics[width=5cm,angle=0]{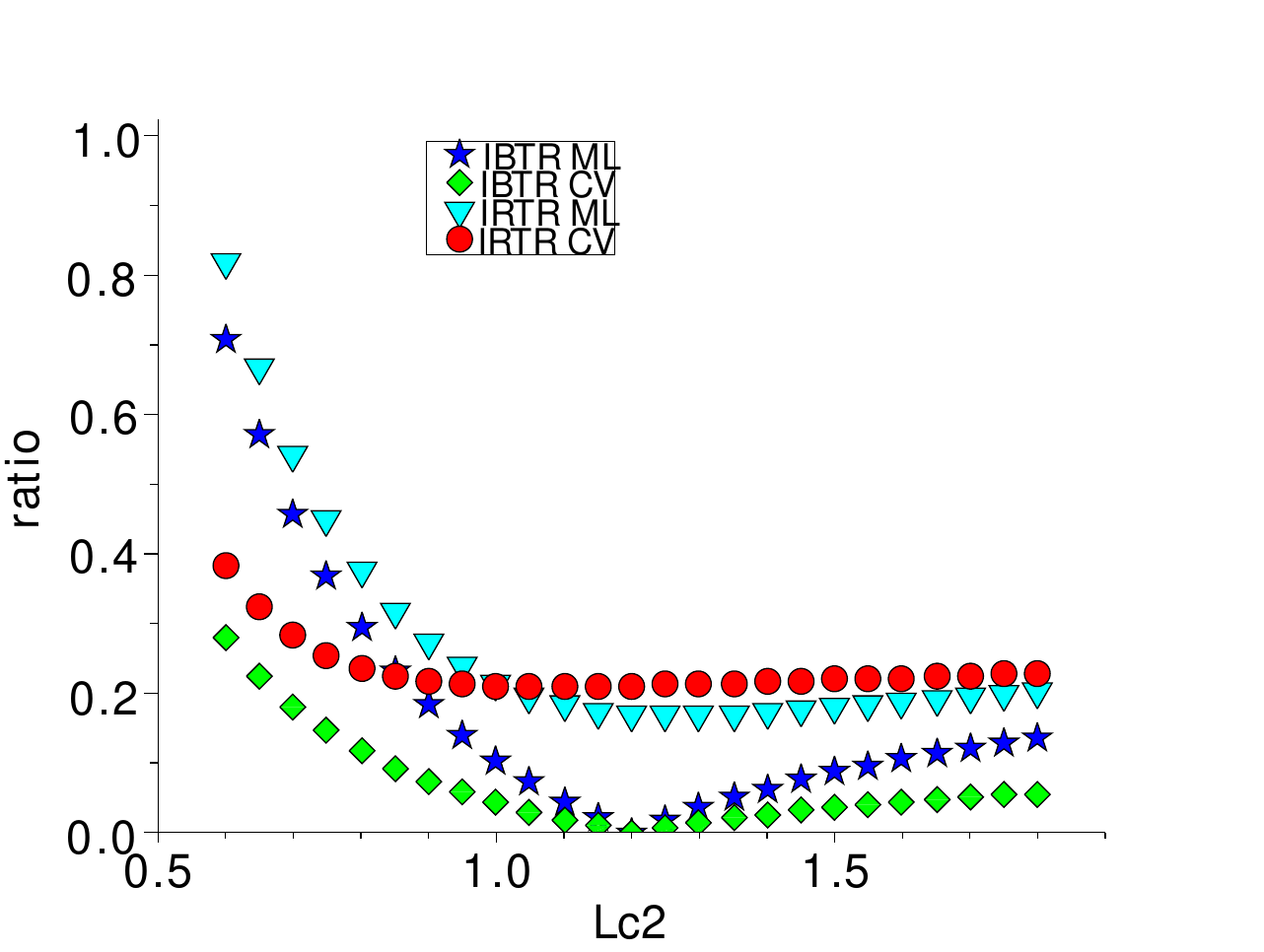}
\end{tabular}
\caption{Same setting as in Figure \ref{fig: IntegreSRS}, but with the LHS-Maximin DOE (see Section \ref{subsection: doe}).
ML is optimal when there is no model error while CV is more robust to model misspecifications.}
\label{fig: IntegreLHS}
\end{figure}

On Figure \ref{fig: IntegreLHS} we plot the results for the LHS-Maximin DOE. The results are similar to these of the SRS DOE.
They also appear to be slightly more pronounced, the IRTR of CV being smaller than the IRTR of ML for a smaller model error.

\begin{figure}[]
\centering
 \hspace*{-2cm}
\begin{tabular}{c c c}
\includegraphics[width=5cm,angle=0]{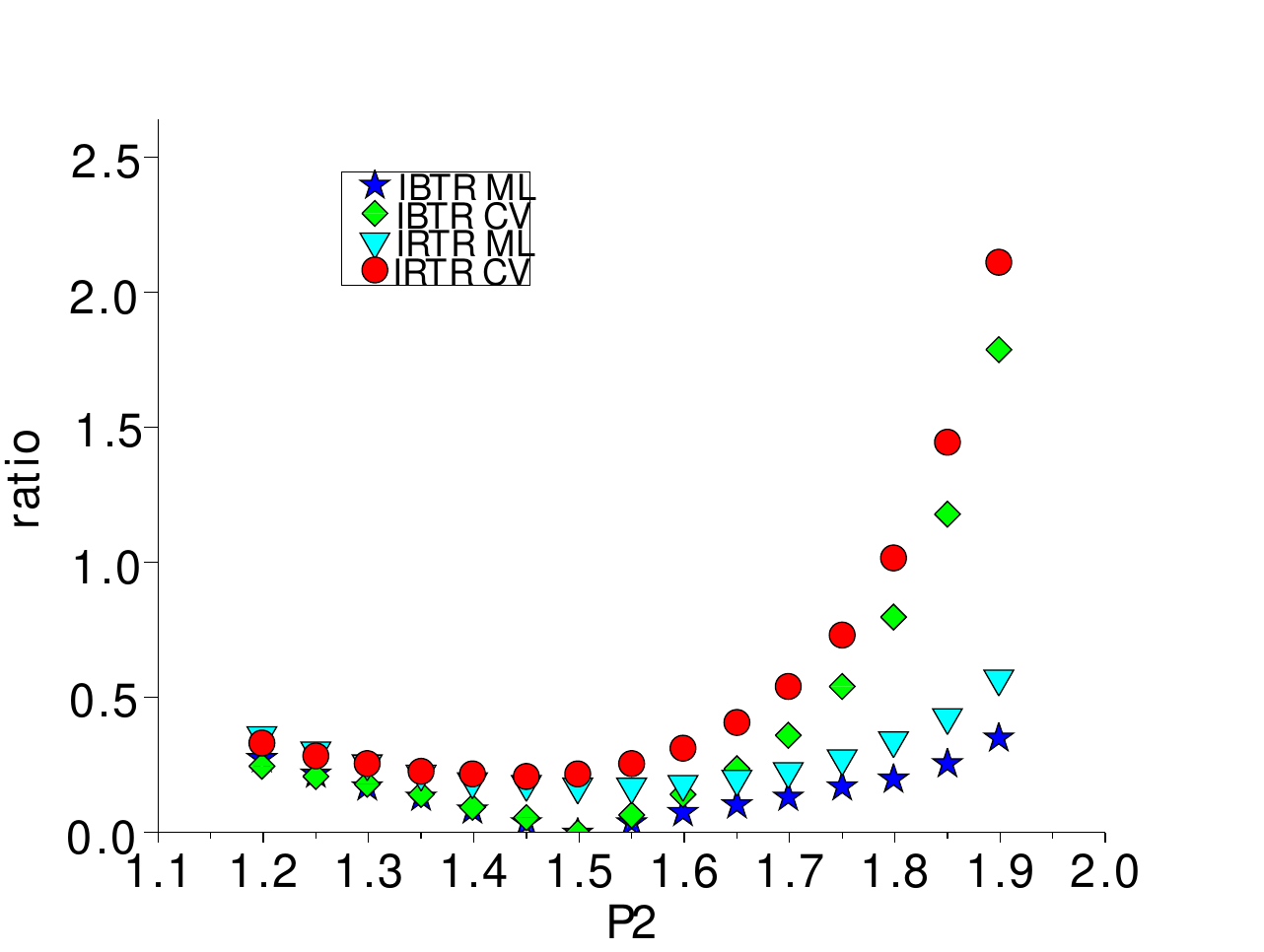} & \includegraphics[width=5cm,angle=0]{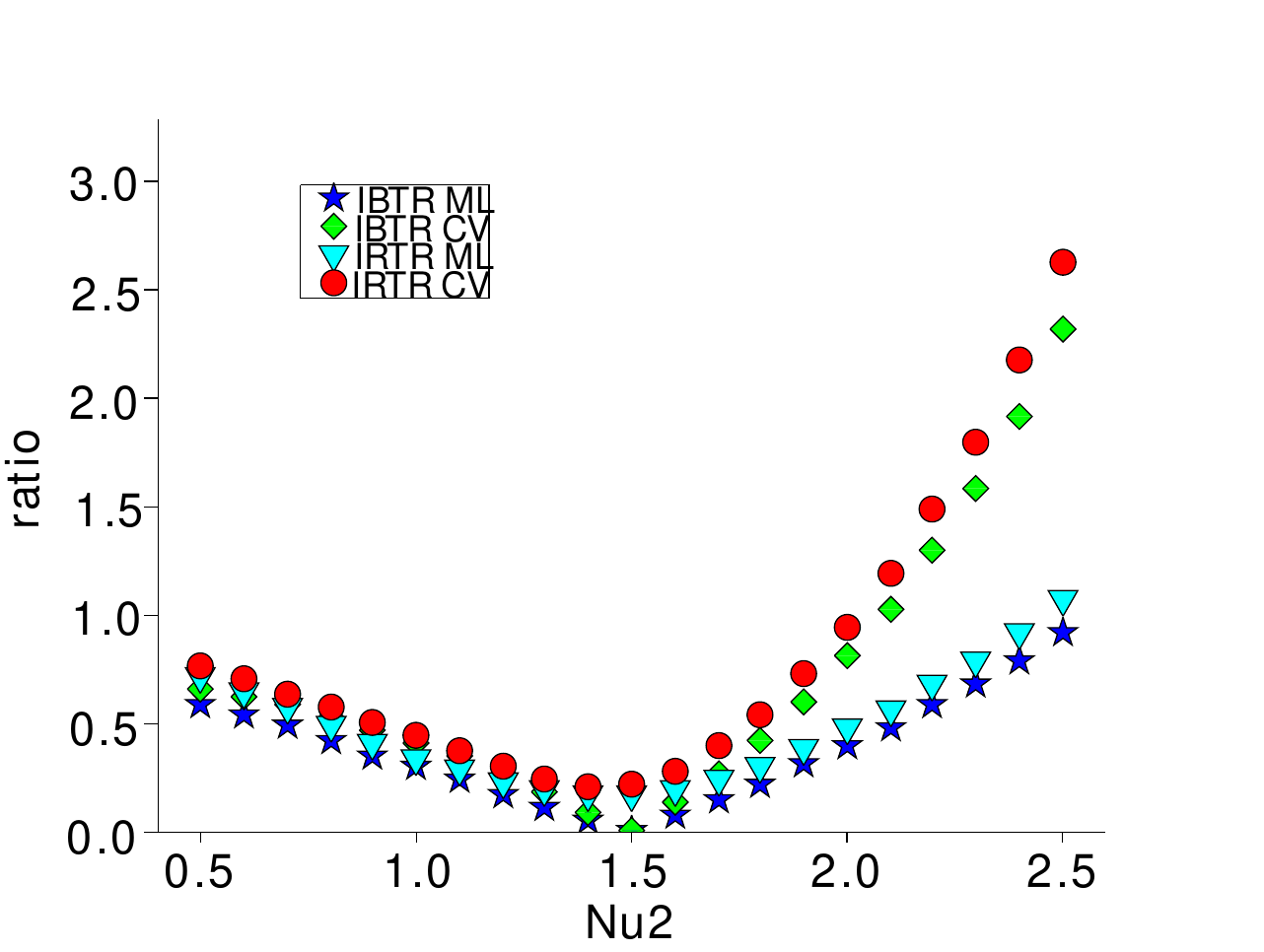}
& \includegraphics[width=5cm,angle=0]{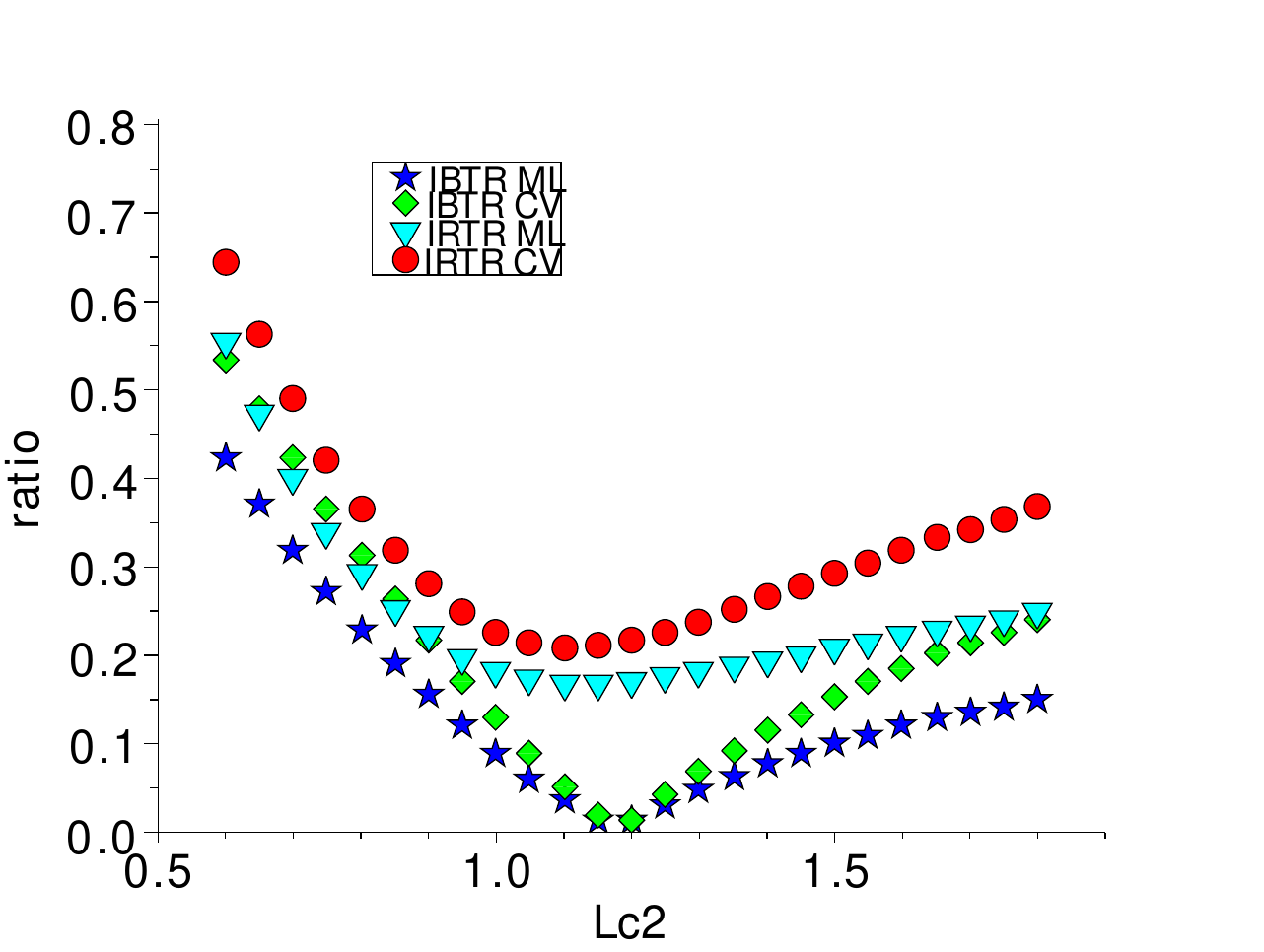}
\end{tabular}
\caption{Same setting as in Figure \ref{fig: IntegreSRS} but with the regular sparse grid DOE (see Section \ref{subsection: doe}).
The results are radically different from the ones obtained with the SRS and LHS-Maximin DOEs. This time CV is less robust to misspecifications of the correlation function.}
\label{fig: IntegreGrid}
\end{figure}

On Figure \ref{fig: IntegreGrid}, we plot the results for the regular grid DOE. The results are radically different from the ones obtained with the SRS and LHS-Maximin
designs. The first comment is that the assessment of predictive variances is much more difficult in case of model misspecification (compare the min, between ML and CV, of IRTR
for the SRS and LHS-Maximin designs against the regular grid). This is especially true for misspecifications on the exponent for the power-exponential correlation function
and on the regularity parameter for the Mat\'ern function. The second comment is that this time CV appears to be less robust than ML to model misspecification.
In particular,
its bias increases faster than ML bias with model misspecification and can be very large. Indeed, having observation points that are on a regular grid, CV estimates
a $\sigma^2$ hyper-parameter adapted only to predictions on the regular grid. Because of the correlation function misspecification, this does not generalize at all to
predictions outside the regular grid. Hence, CV is efficient to assess predictive variances at the points of the regular grid but not to assess predictive variances outside
the regular grid. This is less accentuated for ML because ML estimates a general-purpose $\sigma^2$ and not a $\sigma^2$ for the purpose of assessing predictive variances
at particular points. Furthermore, it is noted in \citet{NSMFVP} that removing a point from a highly structured DOE breaks its structure, which may yield
overpessimistic CV results.

We conclude from these numerical results that, for the SRS and LHS-Maximin designs of
experiments, CV is more robust to model misspecification. It is the contrary for the regular grid, for the structural reasons presented above. This being
said, we do not consider the regular grid anymore in the following numerical results and only consider the SRS and LHS-Maximin designs. Let us finally notice
that the regular grid is not particularly a Kriging-oriented DOE. Indeed, for instance, for $n=71$, it
remains only $17$ distinct points when projecting on the first two base vectors.

\subsection{Influence of the number of points} \label{subsection: influence_n}

Using the same procedure as in Section \ref{subsection: inf_erreur_modele}, we still set $d=5$ and we vary the learning sample size $n$. The pair $R_1,R_2$ is fixed in
the three following different cases. First, $R_1$ is power-exponential $((1.2,...,1.2),1.5)$ and
$R_2$ is power-exponential $((1.2,...,1.2),1.7)$.
Second, $R_1$ is Mat\'ern $((1.2,...,1.2),1.5)$ and $R_2$ is Mat\'ern $((1.2,...,1.2),1.8)$.
Finally, $R_1$ is Mat\'ern $((1.2,...,1.2),1.5)$ and $R_2$ is Mat\'ern $((1.8,...,1.8),1.5)$.
This time, we do not consider integrated quantities of interest and focus on the prediction on the point $x_0$ having all
its components set to $\frac{1}{2}$ (center of domain). 

\begin{figure}[]
\centering
 \hspace*{-2cm}
\begin{tabular}{c c c}
\includegraphics[width=5cm,angle=0]{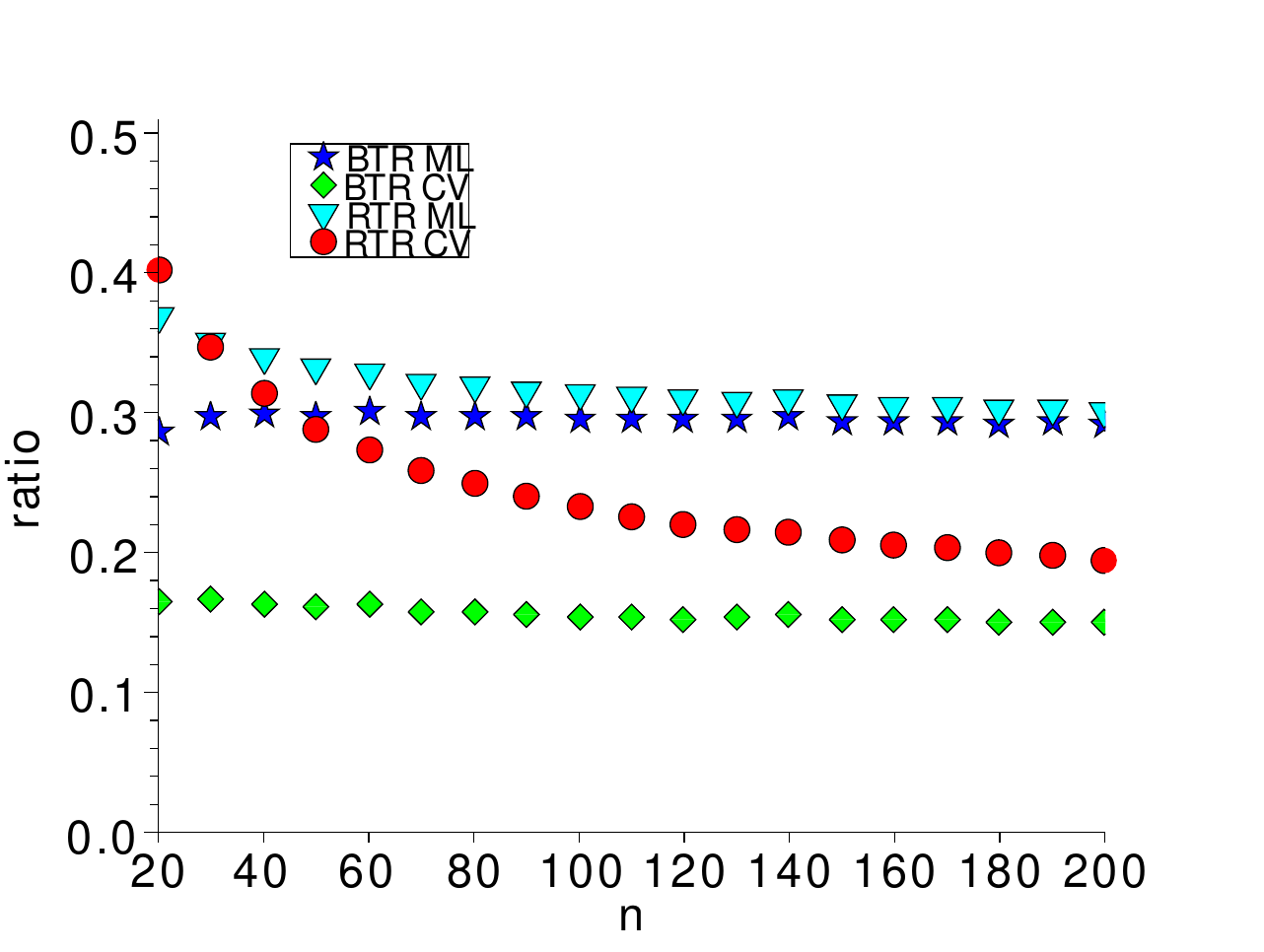} & \includegraphics[width=5cm,angle=0]{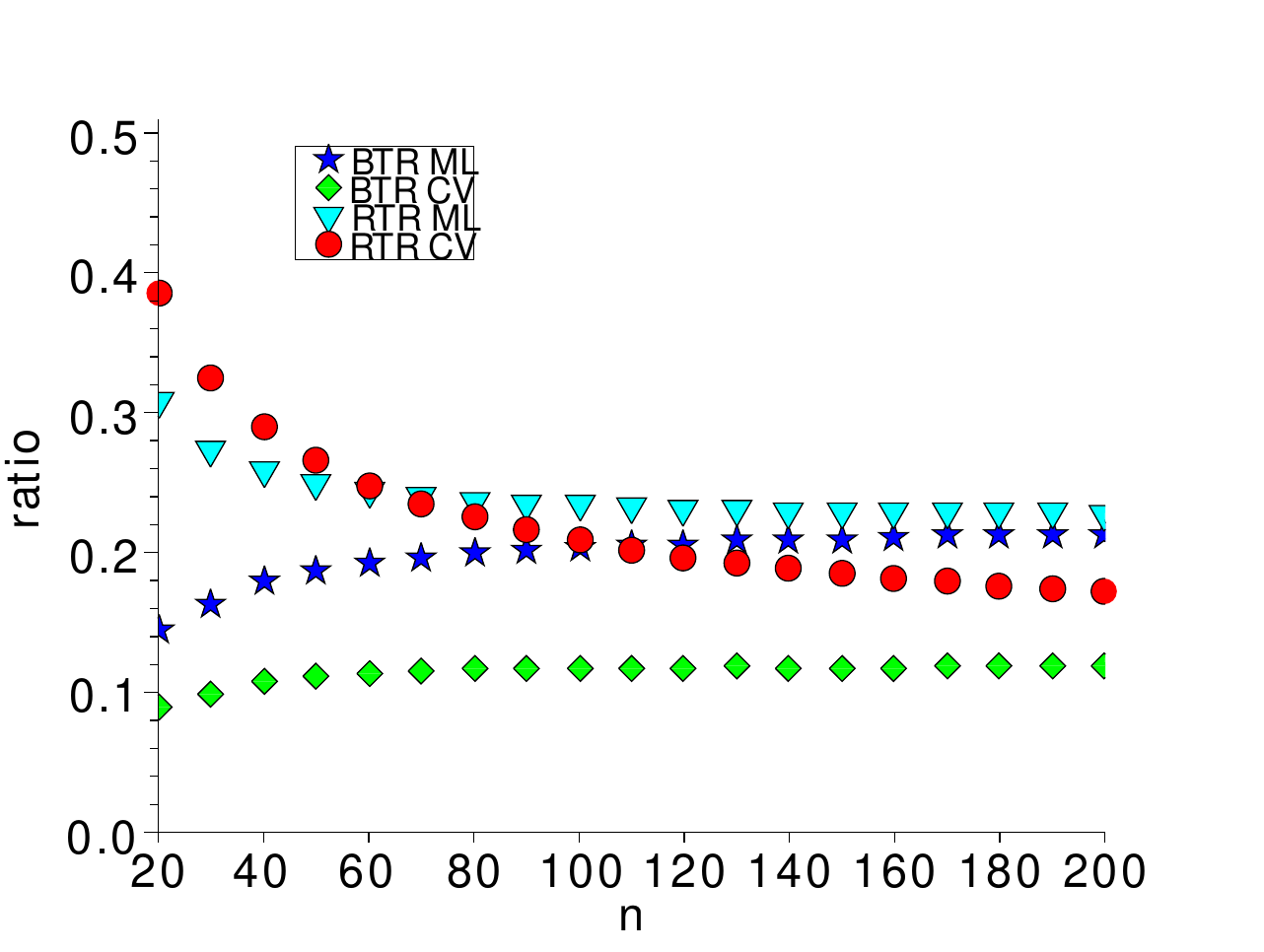}
& \includegraphics[width=5cm,angle=0]{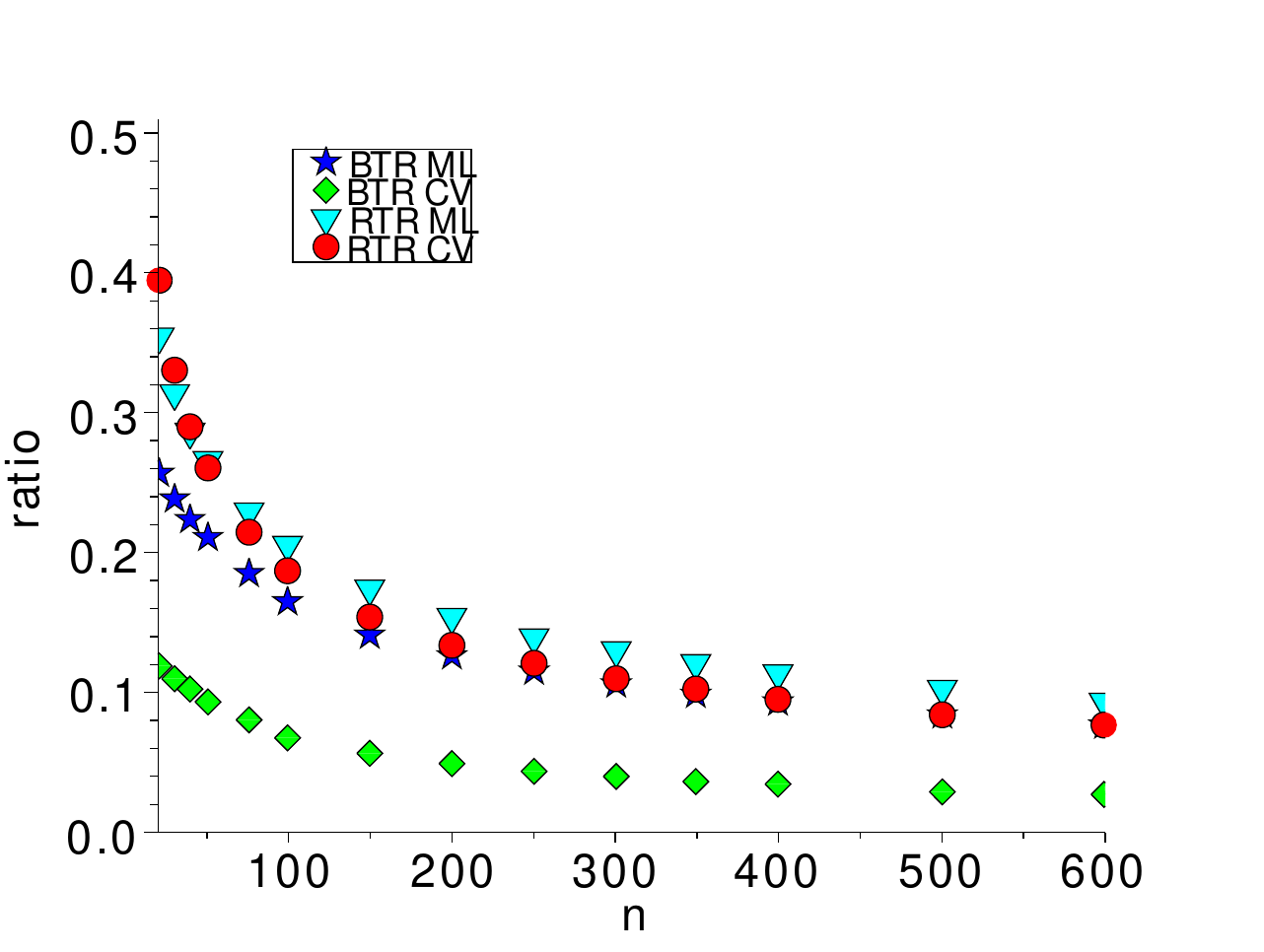}
\end{tabular}
\caption{Influence of the number $n$ of observation points for the SRS DOE (see Section \ref{subsection: doe}). Plot of the RTR and BTR criteria (see Section
\ref{subsection: criteria}) for prediction at the center of the domain and for ML and CV. Left: power-exponential correlation function with error on the
exponent, the true exponent is $p_1=1.5$ and the model exponent
is $p_2 = 1.7$. Middle: Mat\'ern correlation function with error on the regularity parameter, the true regularity parameter is $\nu_1=1.5$ and the model
regularity parameter is $\nu_2=1.8$. Right: Mat\'ern correlation function ($\nu = \frac{3}{2}$) with error on the correlation length, the true correlation length is $\ell_1=1.2$ and the
model correlation length is $\ell_2 = 1.8$.}
\label{fig: InfluenceNSRS}
\end{figure}

On Figure \ref{fig: InfluenceNSRS} we plot the results for the SRS DOE. The first comment is that, as $n$ increases, the BTR does not vanish,
but seems to reach a limit value. This limit value is smaller for CV for the three pairs $R_1,R_2$. Recalling from \eqref{eq: link_RTR_BTR}
that RTR is the sum of BTR and of a relative variance term, we observe that this relative
variance term decreases and seems to vanish when $n$ increases (because BTR becomes closer to RTR). The decrease is much slower for the error on the correlation length
than for the two other errors on the
correlation function. Furthermore, the relative variance term decreases more slowly for CV than for ML. Finally, because CV is better
than ML for the BTR and worse than ML for the relative variance, and because the contribution of BTR to RTR increases with $n$, the ratio of the RTR of ML over the
RTR of CV increases with $n$. This ratio can be smaller than $1$ for very small $n$ and eventually becomes larger than $1$ as $n$ increases
(meaning that CV does better than ML).

\begin{figure}[]
\centering
 \hspace*{-2cm}
\begin{tabular}{c c c}
\includegraphics[width=5cm,angle=0]{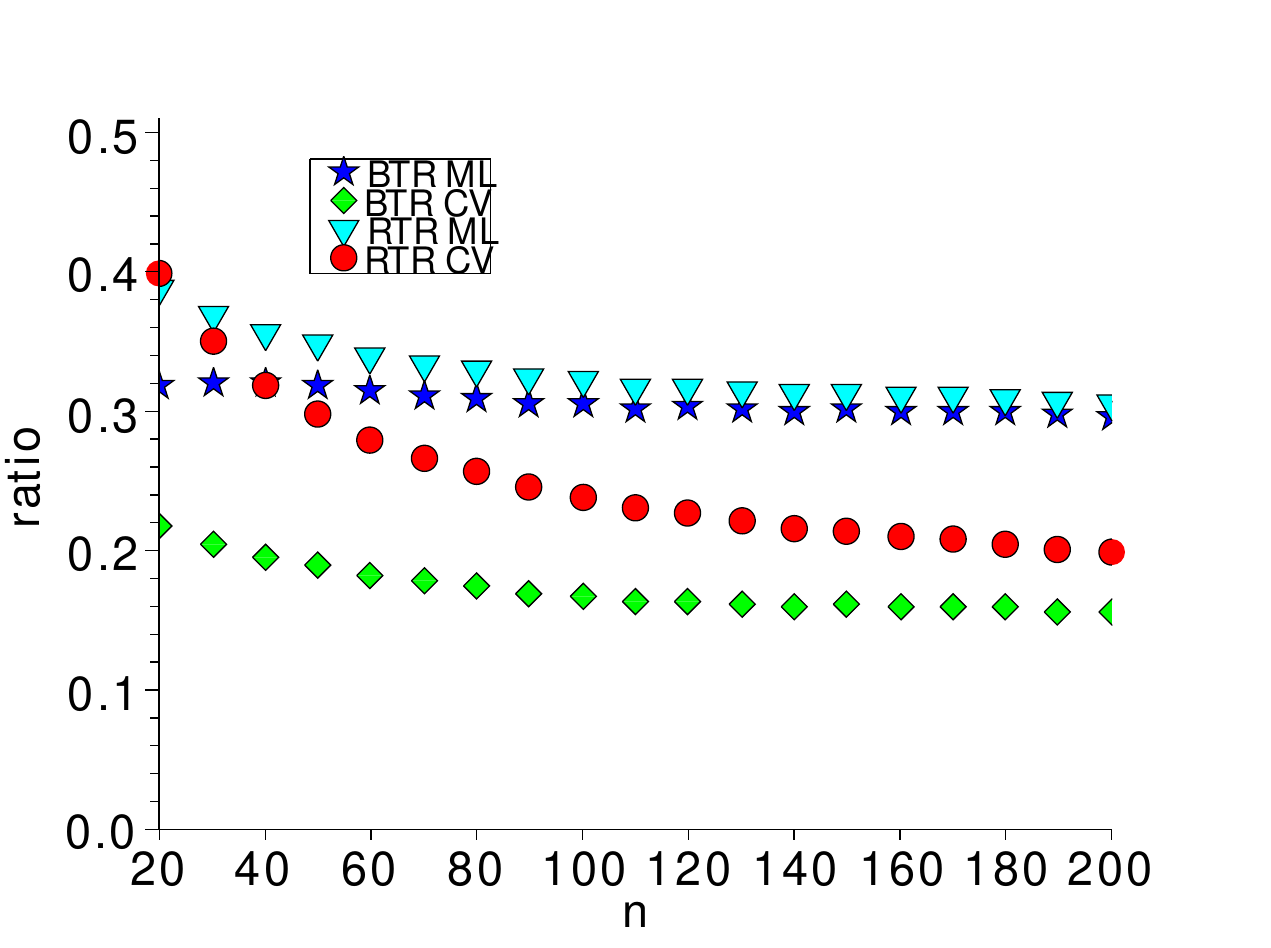} & \includegraphics[width=5cm,angle=0]{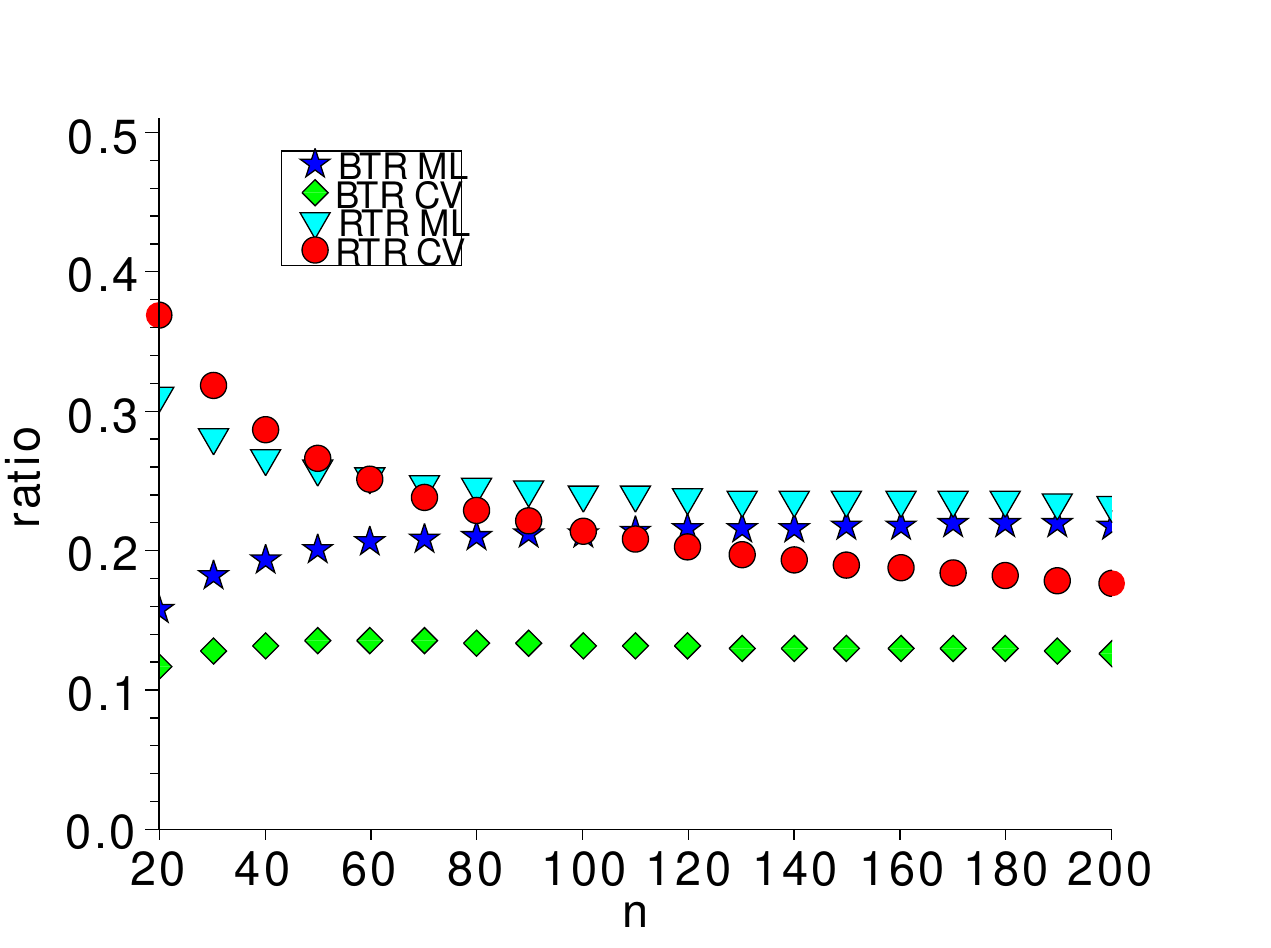}
& \includegraphics[width=5cm,angle=0]{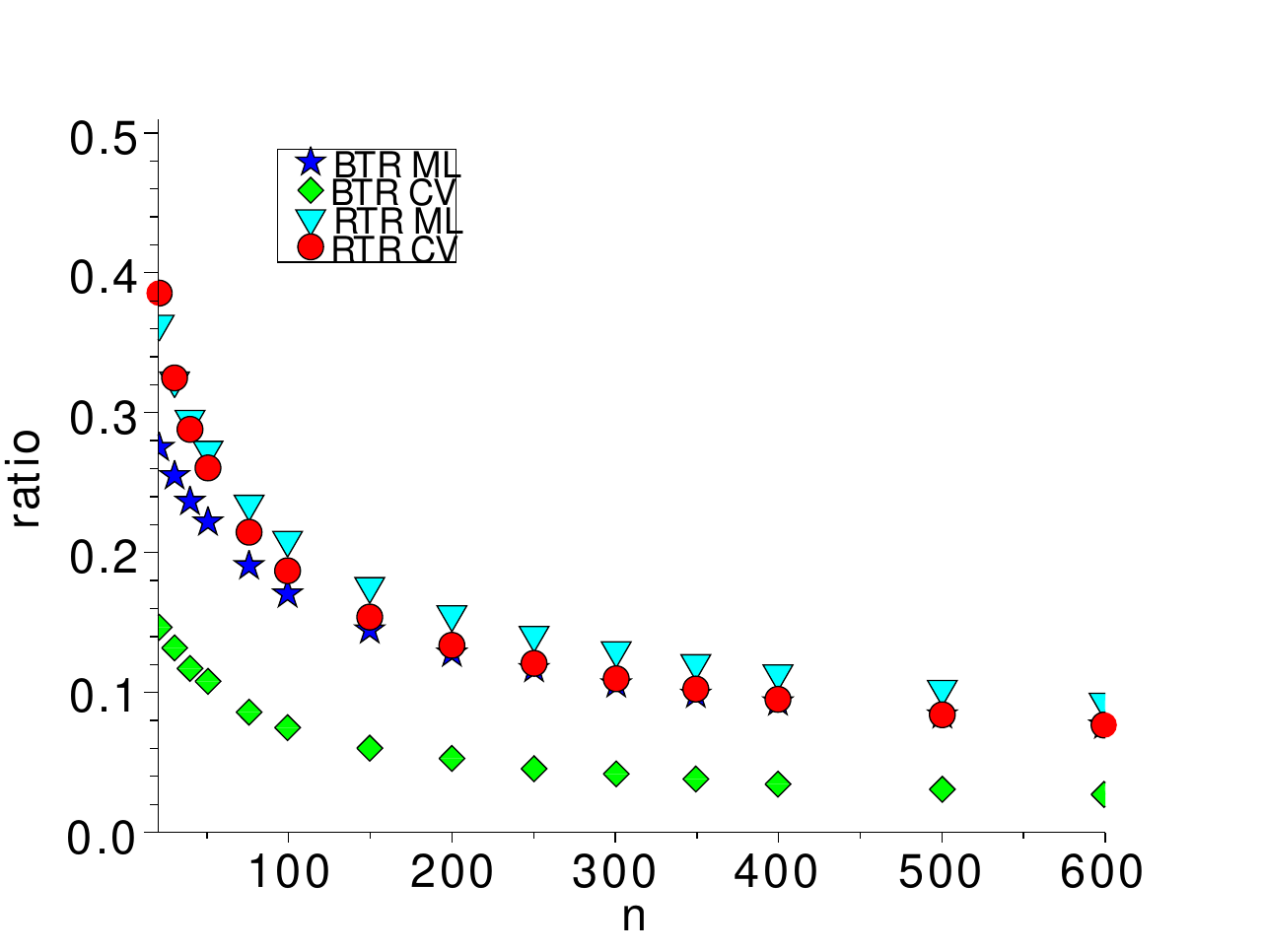}
\end{tabular}
\caption{Same setting as in Figure \ref{fig: InfluenceNSRS}, but with the LHS-Maximin DOE (see Section \ref{subsection: doe}).}
\label{fig: InfluenceNLHS}
\end{figure}

On Figure \ref{fig: InfluenceNLHS} we plot the results for the LHS-Maximin DOE. The results are similar to these of the SRS DOE.
The RTR of CV is smaller than the RTR of ML for a slightly smaller $n$. This confirms the results of Section \ref{subsection: inf_erreur_modele} where
the model error for which the IRTR of ML reaches the IRTR of CV is smaller for LHS-Maximin than for SRS.

\FloatBarrier

\section{Study on analytical functions} \label{section: analyticalfunctions}

In this section, we compare ML and CV on analytical functions, instead of realizations of Gaussian processes, as
was the case in Section \ref{section: numerical_results_variance}. The first goal is to illustrate the results of Section \ref{section: numerical_results_variance}
on the estimation of the variance
hyper-parameter. Indeed, the study of Section \ref{section: numerical_results_variance} is more related to the theory of Kriging (we work on Gaussian processes) while
this section is more related to the application of Kriging (modeling of deterministic functions as realizations of Gaussian processes).
The second goal is to generalize Section
\ref{section: numerical_results_variance} to the case where correlation hyper-parameters are estimated from data. 

\subsection{ML and CV estimations of covariance hyper-parameters} \label{subsection: estimation_correlation}

We consider a set of observations $(x_{1},y_{1}),...,(x_{n},y_{n})$ as in Section \ref{section: framework_variance}, and the family
$\left\{ \sigma^2 R_{\theta}, \sigma^2 >0 , \theta \in \Theta \right\}$ of stationary
covariance functions, with $R_{\theta}$ a stationary correlation function, and $\Theta$ a finite-dimensional set.
We denote by $\EE_{\theta}$ and $var_{\theta}$ the means and variances
with respect to the distribution of a stationary Gaussian process with mean zero, variance one and correlation
function $R_{\theta}$.

The ML estimate of $(\sigma^2,\theta)$ is
$\hat{\theta}_{ML} \in \argmin_{\theta \in \Theta} | \mathbf{\Gamma}_{\theta} |^{1/n} \hat{\sigma}_{ML}^2 ( \mathbf{\Gamma}_{\theta} )$ (see e.g \citet{EMMCCCGP}), with
$\mathbf{\Gamma}_{\theta}$ the correlation matrix of the training sample with correlation function $R_{\theta}$, and
$\hat{\sigma}_{ML}^2 ( \mathbf{\Gamma}_{\theta} )$ as in \eqref{eq: hatsigmaML}. $\hat{\sigma}_{ML}^2 ( \mathbf{\Gamma}_{\hat{\theta}_{ML}}  )$ is the estimate of $\sigma^2$.

For CV we choose to use the following estimator for the hyper-parameter $\theta$:
\begin{equation}  \label{eq: thetaCV}
\hat{\theta}_{CV} \in \argmin_{\theta \in \Theta} \sum_{i=1}^n ( y_{i} - \hat{y}_{i,\theta} (y_{-i}) )^2,
\end{equation} 
with $\hat{y}_{i,\theta} (y_{-i})  = \EE_{\theta} ( y_i | y_1,...,y_{i-1},y_{i+1},...,y_n )$. This is the Leave-One-Out mean square error
criterion that is used, for instance, in \citet{TDACE} when the ML and CV estimations of $\theta$ are compared. Given $\hat{\theta}_{CV}$ and hence
$\mathbf{\Gamma}_{\hat{\theta}_{CV}}$, $\sigma^2$ is estimated by \eqref{eq: hatsigmaCV}.

After elementary transformations for the estimator $\hat{\theta}_{ML}$ and the use of Proposition \ref{prop: dubrule} for
the estimator $\hat{\theta}_{CV}$, the functions to minimize are
\begin{equation} \label{eq: fML}
f_{ML}(\theta) = \frac{1}{n} \log{ det (\mathbf{\Gamma}_{\theta} )}  + \log{ \left( y^t \mathbf{\Gamma}_{\theta}^{-1} y \right) }
\end{equation}
and
\begin{equation} \label{eq: fCV}
f_{CV}(\theta) = y^t \mathbf{\Gamma}_{\theta}^{-1} diag( \mathbf{\Gamma}_{\theta}^{-1} )^{-2} \mathbf{\Gamma}_{\theta}^{-1} y.
\end{equation}
In the supplementary material, we recall
the closed-form expressions of the gradients of $f_{ML}$ and $f_{CV}$,
as functions of the first-order derivatives of the
correlation function. These expressions are available in the literature, see e.g \citet{MLEMRCSR} for ML and \citet{PACHGP} for CV.
The evaluations of the two functions and their gradients have similar computational complexities of the order of $O(n^3)$.

Once we have the closed-form expressions of the gradients at hand, our optimization procedure is based on the Broyden-Fletcher-Goldfarb-Shanno
(BFGS) quasi-Newton
optimization method, implemented in the Scilab function \textit{optim}. Since the functions $f_{ML}$ and $f_{CV}$ may have multiple local minima, the BFGS
method is run several times, by taking the initial points in a LHS-Maximin design. The presence of multiple local minima is discussed e.g in \citet{UKMADCM}. An important point
is that, when $\theta$ is a correlation length, we recommend to use its logarithm to run the optimization. Indeed a correlation length acts as a multiplier in
the correlation, so that using its log ensures that a given perturbation has the same importance, whether applied to a large or a small correlation length. Furthermore, when
one wants to explore uniformly the space of correlation lengths, as is the case with a LHS design, using directly the correlation lengths may give too much emphasis on
large correlation lengths, which is avoided by using their log. 

Another important issue is the numerical inversion of the correlation matrix.
This issue is even more significant when the correlation matrix is ill-conditioned, which happens when the correlation function is smooth (Gaussian or Mat\'ern with a large
regularity parameter). To tackle this issue we recommend
to use the Nugget effect. More specifically, for a given correlation matrix $\mathbf{\Gamma}_{\theta}$, we actually compute
$\mathbf{\Gamma}_{\theta} + \tau^2 I_n$,
with $\tau^2 = 10^{-8}$ in our simulations. A detailed analysis of the influence of the Nugget effect on the hyper-parameter estimation and on the Kriging prediction
is carried out in \citet{ENGPECM}.
However, for the CV estimation of $\theta$, when the correlation function belongs to the Gaussian family, or the Mat\'ern family with large regularity parameter, another structural
problem appears. 
For $\hat{\sigma}^2_{CV}$ very large, as the overall predictive variance term $ \hat{\sigma}^2_{CV} ( 1 - \gamma_{\theta} \mathbf{\Gamma}_{\theta}^{-1} \gamma_{\theta}  ) $
has the same order of magnitude as the observations, the term $ 1 - \gamma_{\theta} \mathbf{\Gamma}_{\theta}^{-1} \gamma_{\theta}  $ is very small.
Hence, a fixed
numerical error on the inversion of $\mathbf{\Gamma}_{\theta}$, however small it is, may cause the term $ 1 - \gamma_{\theta} \mathbf{\Gamma}_{\theta}^{-1} \gamma_{\theta}  $
to be negative. This is what we observe
for the CV case when fitting e.g the correlation lengths of a Gaussian correlation function. The heuristic scheme is that large correlation lengths are estimated, which yields
large $\hat{\sigma}^2_{CV}$, which yields small $(1 - \gamma_{\theta} \mathbf{\Gamma}_{\theta}^{-1} \gamma_{\theta})$, so possibly negative ones.
Notice however that the relative errors of the Kriging prediction terms
$\gamma_{\theta}^t \mathbf{\Gamma}_{\theta}^{-1} y$ are correct. It is noted in \citet[p.7]{UKMADCM} that CV may overestimate
correlation lengths. Hence, to have appropriate predictive variances, one has to ensure that the estimated correlation lengths
are not
too large. Two possible solutions are to penalize either too large correlation lengths or too large $\hat{\sigma}^2_{CV}$ in the minimization of $f_{CV}$.
We choose here the second solution because our experience is that the ideal penalty on the correlation lengths, both ensuring reliable predictive variance computation and having
a minimal effect on the $\theta$ estimation, depends on the DOE substantially.
In practice, we use a penalty for $\hat{\sigma}^2_{CV}$ starting at $1000$ times the empirical variance $\frac{1}{n} y^ty $. This penalty is needed
only for CV when the correlation function is Gaussian or Mat\'ern with free regularity parameter.

\subsection{Procedure}

We consider a deterministic function $f$ on $[0,1]^d$. We generate $n_p=100$ LHS-Maximin training samples of the form
$x_{a,1},f(x_{a,1}),...,x_{a,n},f(x_{a,n})$. We denote $y_{a,i} = f(x_{a,i})$. From each training sample, we estimate $\sigma^2$ and $\theta$
with the ML and CV methods presented in Section \ref{subsection: estimation_correlation}.

We are interested in two criteria on a large Monte Carlo test sample
$x_{t,1},f(x_{t,1}),...,x_{t,n_{t}},f(x_{t,n_{t}})$ on $[0,1]^d$ ($n_t = 10000$).
We denote $y_{t,i} = f(x_{t,i})$,
$\hat{y}_{t,i}(y_a) = \EE_{\hat{\theta}} ( y_{t,i} | y_a )$ and $\hat{\sigma}^2 c_{t,i}^2 (y_a) = \sigma^2 var_{\hat{\theta}} ( y_{t,i} | y_a )$,
where $\hat{\theta}$ comes from either the ML or CV method.

The first criterion is the Mean Square Error (MSE) and evaluates the prediction capability of the estimated correlation function $R_{\hat{\theta}}$:
\begin{equation} \label{eq: MSECrit}
\frac{1}{n_t} \sum_{i=1}^{n_t} ( y_{t,i} - \hat{y}_{t,i} (y_{a}) )^2.
\end{equation}

The second criterion is the Predictive Variance Adequation (PVA):
\begin{equation} \label{eq: PVA}
\left| \log \left(  \frac{1}{n_t} \sum_{i=1}^{n_t} \frac{( y_{t,i} - \hat{y}_{t,i} (y_{a}) )^2}{ \hat{\sigma}^2 c_{t,i}^2 (y_a)} \right) \right|.
\end{equation}
This criterion evaluates the quality of the predictive variances given by the estimated covariance hyper-parameters $\hat{\sigma}^2,\hat{\theta}$.
The smaller the PVA is, the better it is because the predictive variances are globally of the same order than the prediction errors, so that
the confidence intervals are reliable. We use
the logarithm in order to give the same weight to relative overestimation and to relative underestimation of the prediction errors.

We finally average the two
criteria over the $n_p$ training samples.

\subsection{Analytical functions studied}

We study the two following analytical functions. The first one, for $d=3$, is the Ishigami function:
\begin{equation} \label{eq: Ishigami}
f(x_1,x_2,x_3) = \sin(- \pi + 2 \pi x_1) + 7 \sin( (-\pi + 2 \pi x_2) )^2 + 0.1 (- \pi + 2 \pi x_3)^4 \sin(- \pi + 2  \pi x_1).
\end{equation}

The second one, for $d=10$, is a simplified version of the Morris function \citep{FSPPCE},
\begin{eqnarray*}
& & f(x)  =  \sum_{i=1}^{10} w_i(x) + \sum_{1 \leq i < j \leq 6} w_i(x) w_j(x) + \sum_{1 \leq i < j < k  \leq 5} w_i(x) w_j(x) w_k(x)\\
&  & \;\;\;\;\;\;\;\;\;\;\;\;\; +  w_1(x) w_2(x) w_3(x) w_4(x), \\
\mbox{with} & & w_i (x) = 
\begin{cases} 
2 \left( \frac{1.1 x_i}{x_i + 0.1} - 0.5 \right),  & \mbox{if } i = 3,5,7 \\
2 (x_i - 0.5) & \mbox{otherwise}.
\end{cases} \\
\end{eqnarray*}

\subsection{Results with enforced correlation lengths} \label{subsection: result_enforced_l}

We work with the Ishigami function, with $n=100$ observation points. For the correlation function family, we study
the exponential and Gaussian families (power-exponential family of \eqref{eq: pwrexp} with enforced $p=1$ for exponential and $p=2$ for Gaussian). 

For each of these two correlation models, we enforce three vectors $\ell$ of correlation lengths for $R$: an arbitrary isotropic correlation length, a well-chosen isotropic correlation
length and three well-chosen correlation lengths along the three dimensions.
To obtain a well-chosen isotropic correlation length, we generate $n_p=100$ LHS-Maximin DOEs, for
which we estimate the
correlation length by ML and CV as in Section \ref{subsection: estimation_correlation}. We calculate each time the MSE on a test sample of size $10000$ and
the well-chosen correlation length is the one with the smallest
MSE among the $2 n_p$ estimated correlation lengths.
The three well-chosen correlation lengths are obtained similarly.
The three vectors of correlation lengths yield an increasing prediction quality.

The results are presented in Table \ref{table: Case1}. The first comment is that, comparing line $1,2$ and $3$ against line $4,5$ and $6$,
the Gaussian family is more appropriate than the exponential one for the Ishigami function.
Indeed, it yields the smallest MSE among the cases when one uses three different correlation lengths, and the PVA is quite small as well.
This could be anticipated since the Ishigami function is smooth, so a Gaussian correlation model (smooth trajectories) is more adapted than an exponential one
(rough trajectories). The second comment is that the benefit obtained by replacing an isotropic correlation length by different correlation lengths is smaller for the exponential
class than for the Gaussian one. Finally, we see that CV yields much smaller PVAs than ML in line $1$, $2$, $3$ and $4$, in the cases when the correlation function is not
appropriate.
For line $6$, which is the most appropriate correlation
function, ML yields a PVA comparable to CV and for line $5$, ML PVA is smaller than CV PVA. All these comments are in agreement with the main result of Section
\ref{section: numerical_results_variance}: The ML estimation of $\sigma^2$ is more appropriate when the correlation function is well-specified while
the CV estimation is more appropriate when the correlation function is misspecified.

\begin{table} 
\begin{center} 
\begin{tabular}{  c  c  c  c  }
\hline
      Correlation model  &     Enforced hyper-parameters  &   MSE     &     PVA    \\                

\hline

    exponential  &    $[ 1, 1, 1]$     &    $2.01$    &      $ML: 0.50$ $CV: 0.20$  \\                 

     exponential   &  $[ 1.3, 1.3, 1.3]$     & $1.94$    &      $ML: 0.46$ $CV: 0.23$    \\                 

 exponential &     $[ 1.20, 5.03, 2.60]$     &     $1.70$  &        $ML: 0.54$ $CV: 0.19$  \\               

 Gaussian   &   $[ 0.5, 0.5, 0.5]$    &  $4.19$    &     $ ML: 0.98$ $CV: 0.35$      \\             

     Gaussian  &   $[ 0.31, 0.31, 0.31]$   &   $2.03$     &     $ML: 0.16$ $CV: 0.23$      \\             

  Gaussian   &    $[ 0.38, 0.32, 0.42]$     &   $1.32 $   &    $  ML: 0.28$ $CV: 0.29 $  \\           
\hline
\end{tabular}
\end{center} 
\caption{ Mean of the MSE and PVA criteria for the Ishigami function for different fixed correlation models. The MSE is the same between ML and CV as the same correlation
function is used. When the correlation model is misspecified, the MSE is large and CV does better than ML for the PVA criterion.}
\label{table: Case1}
\end{table}

\FloatBarrier

\subsection{Results with estimated correlation lengths}

We work with the Ishigami and Morris functions, with $n=100$ observation points. We use the exponential and Gaussian models as in Section
\ref{subsection: result_enforced_l}, as well as the Mat\'ern model of \eqref{eq: Rmat}. We distinguish two subcases for the vector $\ell$ of correlation lengths.
In \textbf{Case i} we estimate a single isotropic correlation length, while in \textbf{Case a} we estimate
$d$ correlation lengths for the d dimensions.

In Table \ref{table: Case2n100}, we present the results.
For both the Ishigami and Morris functions, the Gaussian model yields smaller MSEs than the exponential model. Indeed, both functions are smooth.
Over the different DOEs, we observe that the estimated Mat\'ern regularity hyper-parameters are large, so that the MSEs and the PVAs for
the Mat\'ern model are similar to the ones of the Gaussian model. Let us notice that for the Ishigami function, the relatively large number $n=100$ of observation
points is required
for the Gaussian correlation model to be more adapted than the exponential one. Indeed, in Table \ref{table: Case2n70}, we show the same results with $n=70$ where the
Gaussian model yields relatively larger MSEs and substantially larger PVAs. Our interpretation is that the linear interpolation that the exponential correlation
function yields can be sufficient, even for a smooth function, if there are not enough observation points.
We also notice that, generally, estimating different correlation lengths (Case a)
yields a smaller MSE than estimating one isotropic correlation length (Case i). In our simulations this is always true except for the Ishigami function with the exponential
model. Indeed, we see in Table \ref{table: Case1} that we get a relatively small benefit for the Ishigami function from using different correlation lengths. Here, this benefit
is compensated by an error in the estimation of the $3$ correlation lengths with $n=100$ observation points.
The overall conclusion is that the Gaussian and Mat\'ern
correlation models are more adapted than the exponential one, and that using different correlation lengths is more adapted than using an isotropic one, provided that there
are enough data to estimate these correlation lengths.

In the exponential case, for both Cases i and a, CV always yields a smaller PVA than ML and yields a MSE that is smaller or similar. In Case a, for the Gaussian
and Mat\'ern correlation functions, the most adapted ones, ML always yields MSEs and PVAs smaller than CV or similar. Furthermore, for the Morris function
with Mat\'ern and Gaussian correlation functions, going from Case i to Case a enhances the advantage of ML over CV.

From the discussion above, we conclude that the numerical experiments yield results,
for the deterministic functions considered here, that are in agreement with the conclusion of Section
\ref{section: numerical_results_variance}: ML is optimal for the best adapted correlation models, while CV is more robust to cases of model misspecification.

\begin{table} 
\begin{center} 
\begin{tabular}{  c  c  c  c  }
\hline
Function  &     Correlation model  &   MSE     &     PVA    \\                
\hline
Ishigami          &          exponential Case i   &       ML: 1.99 CV: 1.97 &         ML: 0.35 CV: 0.23  \\          
Ishigami           &         exponential Case a    &      ML: 2.01 CV: 1.77  &        ML: 0.36 CV: 0.24    \\      
Ishigami            &        Gaussian Case i       &   ML: 2.06 CV: 2.11       &   ML: 0.18 CV: 0.22          \\
Ishigami             &       Gaussian Case a       &   ML: 1.50 CV: 1.53  &        ML: 0.53 CV: 0.50  \\        
Ishigami              &      Mat\'ern Case i  &        ML: 2.19 CV: 2.29     &     ML: 0.18 CV: 0.23      \\    
Ishigami               &     Mat\'ern Case a   &       ML: 1.69 CV: 1.67      &    ML: 0.38 CV: 0.41        \\  
Morris        &            exponential Case i &         ML: 3.07 CV: 2.99  &        ML: 0.31 CV: 0.24  \\        
Morris         &           exponential Case a  &        ML: 2.03 CV: 1.99   &       ML: 0.29 CV: 0.21    \\      
Morris          &          Gaussian Case i      &    ML: 1.33 CV: 1.36   &       ML: 0.26 CV: 0.26         \\ 
Morris           &         Gaussian Case a       &   ML: 0.86 CV: 1.21   &       ML: 0.79 CV: 1.56   \\       
Morris            &        Mat\'ern Case i  &        ML: 1.26 CV: 1.28    &      ML: 0.24 CV: 0.25       \\   
Morris             &       Mat\'ern Case a   &       ML: 0.75 CV: 1.06     &     ML: 0.65 CV: 1.43  \\
\hline
\end{tabular}
\end{center} 
\caption{$n=100$ observation points. Mean of the MSE and PVA criteria over $n_p=100$ LHS-Maximin DOEs for the Ishigami ($d=3$) and Morris ($d=10$)
functions for
different fixed correlation models. When the model is misspecified, the MSE is large and the CV does better compared to ML for the MSE and PVA criterion.}
\label{table: Case2n100}
\end{table}

\begin{table} 
\begin{center} 
\begin{tabular}{  c  c  c  c  }
\hline
Function  &     Correlation model  &   MSE     &     PVA    \\                
\hline

Ishigami   &         exponential Case a  &        ML: 3.23 CV: 2.91  &        ML: 0.27 CV: 0.26   \\       
Ishigami      &     Gaussian Case a    &       ML: 3.15 CV: 4.13    &      ML: 0.72 CV: 0.76   \\

\hline
\end{tabular}
\end{center} 
\caption{$n=70$ observation points. Mean of the MSE and PVA criteria over $n_p=100$ LHS-Maximin DOEs for the Ishigami ($d=3$) and Morris ($d=10$) functions for
the exponential correlation model. Contrary to the case $n=100$ of Table \ref{table: Case2n100}, the Gaussian correlation model does not yield smaller MSEs than
the exponential one.}
\label{table: Case2n70}
\end{table}

\subsection{Case of universal Kriging}

So far, the case of simple Kriging has been considered, for which the underlying Gaussian process is considered centered. The case of universal Kriging can equally
be studied, for which
this process is considered to have a mean at location $x$ of the form $\sum_{i=1}^p \beta_i g_i(x)$, with known functions $g_i$ and unknown coefficients $\beta_i$. For instance
a closed-form formula similar to that of Proposition \ref{prop: main} can be obtained in the same fashion, and virtual LOO formulas are also
available \citep{CVKUN}. We have chosen to focus on the simple Kriging case because
we are able to address as precisely as possible the issue of the covariance function class misspecification, the Kriging model depending only on the covariance function
choice. Furthermore it is shown in \citet[p.138]{ISDSTK} that the issue of the mean function choice for the Kriging model is much less crucial than that of the
covariance function choice.

Nevertheless, in Table \ref{table: H} we study, for the Ishigami function, the influence of using a universal Kriging model
with either a constant or affine mean function. The process is the same as for Table \ref{table: Case2n100}. We first see that using a non-zero mean
does not improve significantly the Kriging model.
It is possible to observe a slight improvement only with the exponential covariance structure,
which we can interpret because a smooth mean function makes the Kriging model more adapted
to the smooth Ishigami function.
On the contrary, for the Gaussian covariance structure, the mean function over-parameterizes
the Kriging model and slightly damages its performances. Let us also notice that CV appears to be more sensitive to this over-parameterization, its MSE increasing
with the complexity of the mean function. This can be observed similarly in the numerical experiments in \citet{UKMADCM}.
The second overall conclusion is that the main finding of Section \ref{section: numerical_results_variance} and of
Table \ref{table: Case2n100} is confirmed: CV has smaller MSE and PVA for the misspecified exponential structure, while ML is optimal for the
Gaussian covariance structure which is the most adapted and yields the smallest MSE.

\begin{table} 
\begin{center}
 \hspace*{-1.8cm} 
\begin{tabular}{  c  c  c  c  c  }
\hline
Function  &   Mean function model &   Correlation model  &   MSE     &     PVA    \\                
\hline

Ishigami   &   constant   &     exponential Case a  &        ML: 1.96 CV: 1.74  &        ML: 0.39 CV: 0.24   \\       
Ishigami      &   affine  &  exponential Case a    &       ML: 1.98 CV: 1.75    &      ML: 0.40 CV: 0.24   \\
Ishigami   &   constant   &     Gaussian Case a  &        ML: 1.54 CV: 1.63  &        ML: 0.54 CV: 0.54   \\       
Ishigami      &   affine  &  Gaussian Case a    &       ML: 1.58 CV: 1.78    &      ML: 0.57 CV: 0.57   \\

\hline
\end{tabular}
\end{center} 
\caption{$n=100$ observation points. Mean of the MSE and PVA criteria over $n_p=100$ LHS-Maximin DOEs for the Ishigami ($d=3$) function and the
exponential and Gaussian correlation models. The incorporation of the mean function does not change the conclusions of Table \ref{table: Case2n100}.}
\label{table: H}
\end{table}

\FloatBarrier

\section{Conclusion}

In this paper, we have carried out a detailed analysis of ML and CV for the estimation of the covariance hyper-parameters
of a Gaussian process, with a misspecified parametric family of covariance functions. This analysis has been carried out using a two-step approach.
We have first studied the estimation of a global variance hyper-parameter, for which the correlation function is misspecified.
In this framework, we can control precisely the degree of model
misspecification and we obtain closed-form expressions for the mismatch indices that we have introduced. We conclude from the numerical
study of these formulas that when the model is misspecified, CV performs better than ML.
Second, we have studied the general case when the
correlation hyper-parameters are estimated from data via numerical experiments on analytical functions. We confirm the results of the first step, and generalize
them. 

Because CV is more robust to model misspecification, it is less likely than ML to yield substantially incorrect predictions.
However, ML is
more likely to yield the best predictions, under the condition that the correlation function family is ideally specified. Hence, CV is a suitable method
for problems when one gives priority to robustness over best possible performance. Studying data-driven methods, to choose between ML and CV
according to a specific objective, may motivate further research.
  
We also pointed out in Section \ref{subsection: inf_erreur_modele} that some DOEs that are too regular are not suitable for CV. Investigating this
further and studying numerical criteria
for quantifying the suitability of a DOE for CV
could motivate further research. One could be also interested in CV-oriented
adaptative improvement of a DOE. Finally, we emphasized the need, in some cases, for a penalization of large correlation lengths for CV. We enforced this
penalization via a penalization of too large estimated global variance hyper-parameters.
The motivation here was purely numerical, but there could be a statistical motivation for doing so.
The statistical analysis of penalization methods, similar to the one proposed here, for CV may motivate further research. In the ML case, there exist penalizations
of correlation hyper-parameters that both reduce the small sample variance, and recover the asymptotic distribution of the non-penalized case,
as shown in \citet{ACEUPLGKM}.

\section*{Acknowledgments}
The author would like to thank his advisors Josselin Garnier (Laboratoire de Probabilit\'es et Mod\`eles Al\'eatoires \& Laboratoire
Jacques-Louis Lions, University Paris Diderot) and
Jean-Marc Martinez (French Alternative Energies and Atomic Energy
Commission - Nuclear Energy Division at CEA-Saclay, DEN, DM2S, STMF, LGLS) for their advice and suggestions. The author also thanks
the anonymous reviewers for their helpful suggestions.
The author presented the contents of the manuscript at a workshop of the ReDICE consortium, where he benefited from constructive comments and suggestions.

\newpage

\FloatBarrier

\appendix

\section{Proof of Proposition \ref{prop: main}} \label{appendix: preuveMain}

Using the definition of the Risk, the expression of $\hat{\sigma}^2$, \eqref{eq: riskdepred} and \eqref{eq: cx0}, we get:
\begin{eqnarray*}
\mathcal{R}_{ \hat{\sigma}^2 , x_0 } &  = & \EE_1 \big{[}  ( \gamma_1^t \mathbf{\Gamma}_1^{-1} y - \gamma_2^t \mathbf{\Gamma}_2^{-1} y )^2 + 1 - \gamma_1^t \mathbf{\Gamma}_1^{-1} \gamma_1  \\
&  & - y^t \mathbf{M} y  (1 - \gamma_2^t \mathbf{\Gamma}_2^{-1} \gamma_2)   \big{]}^2 \\
& = & \EE_1 \big{[} y^t (\mathbf{\Gamma}_2^{-1} \gamma_2 - \mathbf{\Gamma}_1^{-1} \gamma_1  )( \gamma_2^t \mathbf{\Gamma}_2^{-1} -  \gamma_1^t \mathbf{\Gamma}_1^{-1} ) y \\
&  & + 1 - \gamma_1^t \mathbf{\Gamma}_1^{-1} \gamma_1 - y^t \mathbf{M} y  (1 - \gamma_2^t \mathbf{\Gamma}_2^{-1} \gamma_2) \big{]}^2. \\
\end{eqnarray*}
Then, writing $y = \mathbf{\Gamma}_1^{ \frac{1}{2} } z$ with $z \sim_1 \N(0,\mathbf{I}_n)$, we get:
\begin{equation} \label{eq: hatsigmafinpreuve}
\mathcal{R}_{ \hat{\sigma}^2 , x_0 } = \EE_1 \left( z^t \tilde{\mathbf{M}}_{0} z + c_1 - c_2 z^t \tilde{\mathbf{M}}_1 z   \right)^2,
\end{equation}
with
\[
\tilde{\mathbf{M}}_{0} = \mathbf{\Gamma}_1^{\frac{1}{2}} (\mathbf{\Gamma}_2^{-1} \gamma_2 - \mathbf{\Gamma}_1^{-1} \gamma_1  )( \gamma_2^t \mathbf{\Gamma}_2^{-1} -  \gamma_1^t \mathbf{\Gamma}_1^{-1} ) \mathbf{\Gamma}_1^{\frac{1}{2}}
\]
and
\[
\tilde{\mathbf{M}}_1 = \mathbf{\Gamma}_1^{\frac{1}{2}}  \mathbf{M} \mathbf{\Gamma}_1^{\frac{1}{2}}.
\]
To compute this expression, we use the following lemma. The proof relies only on $4$\textsuperscript{th} moment calculation for Gaussian
variables and is omitted.
\begin{lem} \label{lem: zAzzBz}
Let $z \sim_1 \N(0,\mathbf{I}_n)$, and $\mathbf{A}$ and $\mathbf{B}$ be $n \times n$ real symmetric matrices. Then:
\[
\EE_1( z^t \mathbf{A} z z^t \mathbf{B} z ) = f(\mathbf{A},\mathbf{B}).
\]
\end{lem}
Using the lemma and expanding \eqref{eq: hatsigmafinpreuve} yields
\begin{eqnarray} \label{eq: Mtilde}
\mathcal{R}_{ \hat{\sigma}^2 , x_0 } & = & f(\tilde{\mathbf{M}}_{0},\tilde{\mathbf{M}}_{0}) + 2 c_1 tr(\tilde{\mathbf{M}}_{0}) - 2 c_2 f(\tilde{\mathbf{M}}_{0},\tilde{\mathbf{M}}_1) \\
&  & + c_1^2 - 2 c_1 c_2 tr(\tilde{\mathbf{M}}_1) + c_2^2 f(\tilde{\mathbf{M}}_1,\tilde{\mathbf{M}}_1). \nonumber
\end{eqnarray}

Finally, based on $tr(\mathbf{A}\mathbf{B}) = tr(\mathbf{B}\mathbf{A})$, we can replace $\tilde{\mathbf{M}}_{0}$ and $\tilde{\mathbf{M}}_{1}$ by $\mathbf{M}_0$
and $\mathbf{M}_1$ in \eqref{eq: Mtilde}, which completes the proof.

\section{On the optimality of $\hat{\sigma}_{ML}^2$ when $R_1 = R_2$} \label{subsection: R1=R2}

Here, we consider the case $R_2=R_1$.
We first recall the Cram\'er-Rao inequality, which states that when
$\sigma^2=1$, for an unbiased estimator $\hat{\sigma}^2$ of $\sigma^2$:
\[
var_1( \hat{\sigma}^2 ) \geq \EE_1^{-1} \left\lbrack \left( \frac{\partial}{\partial \sigma^2} \left( \ln{ L(y,\sigma^2) } \right)_{\sigma^2=1}   \right)^2 \right\rbrack,
\]
with, $L(y,\sigma^2) \propto \frac{1}{(\sigma^2)^{\frac{n}{2}}} \exp{ \left(- \frac{y^t \mathbf{\Gamma}_1^{-1} y }{2 \sigma^2} \right) }$, the likelihood of the observations.
We then calculate the Cram\'er-Rao bound:
\begin{eqnarray*}
\EE_1^{-1} \left \lbrack \left( \frac{\partial}{\partial \sigma^2} \left( \ln{ L(y,\sigma^2) } \right)_{\sigma^2=1}   \right)^2 \right \rbrack
& = &  \EE_1^{-1} \left\lbrack \left( \frac{\partial}{\partial \sigma^2} \left( - \frac{n}{2} \ln{\sigma^2} - \frac{ y^t \mathbf{\Gamma}_1^{-1} y }{2 \sigma^2} \right)_{\sigma^2=1}  \right)^2 \right\rbrack  \\
& = & \EE_1^{-1} \left \lbrack  \frac{n^2}{4} + \frac{1}{4} ( y^t \mathbf{\Gamma}_1^{-1} y )^2 - \frac{n}{2} y^t \mathbf{\Gamma}_1^{-1} y \right \rbrack \\
& = & \left( \frac{n^2}{4} + \frac{ n^2 + 2n }{4} - \frac{n^2}{2}  \right)^{-1}  \\
& = & \frac{2}{n},
\end{eqnarray*}
where we used Lemma \ref{lem: zAzzBz} with $\mathbf{A}=\mathbf{B}=\mathbf{I}_n$ to show $\EE_1 \left( ( y^t \mathbf{\Gamma}_1^{-1} y )^2 \right) = n^2 + 2 n$.
Hence, the Cram\'er-Rao bound of the statistical model $\C$ is $\frac{2}{n}$ when $\sigma^2=1$.

Finally,
$var_1(\hat{\sigma}^2_{ML}) = var_1( \frac{1}{n} y^t \mathbf{\Gamma}_1^{-1} y ) = \frac{1}{n^2} var_1( \sum_{i=1}^n z_i^2 ) = \frac{2}{n}$, with
$z = \mathbf{\Gamma}_1^{-\frac{1}{2}} y \sim_1 \N(0,\mathbf{I}_n)$.

\section{Convergence of \eqref{eq: critLOO}} \label{subsection: preuveCritLOO}

\begin{prop} \label{prop: critLOO}
Assume $\X = \RR^d$ and that the observation points are distinct and located on the infinite regular
grid ${ \{ ( i_1 \delta,...,i_d \delta ) , (i_1,...,i_d) \in \Z^d \} }$ for
fixed $\delta > 0$. When $\sigma^2 R_2$ is the covariance function of the Gaussian process $Y$, is summable on the regular grid, and has a strictly positive continuous
Fourier transform, then the criterion $C_{LOO}$ has mean one and converges in mean square to one as $n \to +\infty$.
\end{prop}

\paragraph{Proof}

Consider, for simplicity, $\sigma^2=1$. Using the formulas of Proposition \ref{prop: dubrule} and introducing
$z = \mathbf{\Gamma}_2^{-\frac{1}{2}} y \sim_2 \N(0,\mathbf{I}_n)$ yields:
\begin{equation*}
C_{LOO} = \frac{1}{n} \sum_{i=1}^n \frac{(y_i - \hat{y}_{i,-i} )^2}{ c_{i,-i}^2} =  \frac{1}{n } z^t \mathbf{\Gamma}_2^{-\frac{1}{2}} \left[ diag( \mathbf{\Gamma}_2^{-1} ) \right]^{-1} \mathbf{\Gamma}_2^{-\frac{1}{2}} z,
\end{equation*}
with $diag( \mathbf{\Gamma}_2^{-1} )$ the matrix obtained by setting to $0$ all non-diagonal terms of $\mathbf{\Gamma}_2^{-1}$. Then,
\begin{equation*}
 \EE_2( C_{LOO} ) =    \frac{1}{n } tr \left( \mathbf{\Gamma}_2^{-\frac{1}{2}} \left[ diag( \mathbf{\Gamma}_2^{-1} ) \right]^{-1} \mathbf{\Gamma}_2^{-\frac{1}{2}} \right)  \\
  =  \frac{1}{n } \sum_{i=1}^n \frac{ (\mathbf{\Gamma}_2^{-1})_{i,i} }{ (\mathbf{\Gamma}_2^{-1})_{i,i} }  \\
  =  1.
\end{equation*}
Furthermore, 
\begin{equation} \label{eq: critLOO_var}
var_2( C_{LOO} ) =  \frac{2}{n^2 } tr \left( \mathbf{\Gamma}_2^{-1} \left[ diag( \mathbf{\Gamma}_2^{-1} ) \right]^{-1} \mathbf{\Gamma}_2^{-1} \left[ diag( \mathbf{\Gamma}_2^{-1} ) \right]^{-1}  \right). 
\end{equation}
Then, with $\lambda_{min}(\mathbf{\Gamma})$ and $\lambda_{max}(\mathbf{\Gamma})$ the smallest and largest eigenvalues of a symmetric, strictly positive,
matrix $\mathbf{\Gamma}$,
\begin{equation} \label{eq: critLOO_diag}
\lambda_{max}( \left[ diag( \mathbf{\Gamma}_2^{-1} ) \right]^{-1} ) = \underset{1 \leq i \leq n}{max} \frac{1}{ ( \mathbf{\Gamma}_{2}^{-1} )_{i,i} } \leq \frac{ 1 }{ \lambda_{min} ( \mathbf{\Gamma}_{2}^{-1} ) } = \lambda_{max}(\mathbf{\Gamma}_2).
\end{equation}
Hence, from \eqref{eq: critLOO_var} and \eqref{eq: critLOO_diag},
\[
 \frac{2}{n} \left( \frac{\lambda_{min}(\mathbf{\Gamma}_2)}{\lambda_{max}(\mathbf{\Gamma}_2)} \right)^2 \leq  var_2( C_{LOO} ) \leq \frac{2}{n} \left( \frac{\lambda_{max}(\mathbf{\Gamma}_2)}{\lambda_{min}(\mathbf{\Gamma}_2)} \right)^2.
\]
Hence $0 < \underset{ n }{inf} \lambda_{min} (\mathbf{\Gamma}_2) \leq \underset{ n }{sup}  \lambda_{max} (\mathbf{\Gamma}_2) < + \infty$ implies the proposition.

First, for all $n$, $\lambda_{max} (\mathbf{\Gamma}_2)$ is smaller than the largest $l^1$ norm of the rows of $\mathbf{\Gamma}_2$ and hence is smaller than the absolute sum of $R_2$ over the
infinite regular grid. Hence we have $\underset{ n }{sup}  \lambda_{max} (\mathbf{\Gamma}_2) < + \infty$. 

Then, let $\hat{R}_2$ be the Fourier transform of $R_2$, and let $R_{min} >0$ be its infimum over $[-\frac{\pi}{\delta} , \frac{\pi}{\delta} ]^d$. For all
$n$, $x_1,...,x_n$ in the infinite regular grid, and for all $\epsilon \in \RR^n$:
\begin{eqnarray*}
\epsilon^t \mathbf{\Gamma}_2 \epsilon & = & \sum_{k,l=1}^n  R_2(x_k-x_l) \epsilon_k \epsilon_l 
 =  \sum_{k,l=1}^n \epsilon_k \epsilon_l \int_{\RR^d} e^{ i t (x_k - x_l) } \hat{R}_2(t) dt 
 =  \int_{\RR^d} \hat{R}_2(t) \left\lvert \sum_{k=1}^n \epsilon_k e^{ i t x_k } \right\rvert^2 dt \\
& \geq & R_{min} \int_{ [-\frac{\pi}{\delta} , \frac{\pi}{\delta} ]^d } \left\lvert \sum_{k=1}^n \epsilon_k e^{ i t x_k } \right\rvert^2 dt
 =  R_{min} \int_{ [-\frac{\pi}{\delta} , \frac{\pi}{\delta} ]^d }  \sum_{k,l=1}^n \epsilon_k \epsilon_l e^{ i t (x_k - x_l) } dt.
\end{eqnarray*}
Then, because $\int_{ [-\frac{\pi}{\delta} , \frac{\pi}{\delta} ] } e^{ i t p \delta } dt = 0$ for a non-zero integer $p$, and because the $x_k$, $1 \leq k \leq n$,
are distinct and
located on the regular grid, we have $\int_{ [-\frac{\pi}{\delta} , \frac{\pi}{\delta} ]^d }   e^{ i t (x_k - x_l) } dt = 0$ for $k \neq l$. Hence:
\begin{equation*}
\epsilon^t \mathbf{\Gamma}_2 \epsilon  \geq  R_{min} \int_{ [-\frac{\pi}{\delta} , \frac{\pi}{\delta} ]^d }  \sum_{k=1}^n \epsilon_k^2 dt
 = R_{min} \left( \frac{2 \pi}{\delta} \right)^d || \epsilon ||^2,  \\
\end{equation*}
which in turn implies $0 < \underset{ n }{inf} \lambda_{min} (\mathbf{\Gamma}_2)$ and thus completes the proof.

\newpage

\bibliographystyle{model2-names}
\bibliography{Biblio}

\end{document}